\documentclass[12pt]{iopart}
\usepackage{latexsym}
\usepackage{color}
\usepackage{hyperref}
\usepackage{amssymb}
\usepackage{comment}

\newtheorem{theorem}{Theorem}

\newtheorem{assumption}[theorem]{Assumption}
\newtheorem{proposition}[theorem]{Proposition}

\newtheorem{remark}[theorem]{Remark}

\newcommand{\eM}{\Phi^{\star}}

\def\IN{\mathbb{N}}

\def\IP{\mathbb{P}}
\def\IR{\mathbb{R}}

\def\IE{\mathbb{E}}

\newcommand{\cs}[1]{{\color{darkblue}{#1}}}

\newcommand{\eN}{N_{\tiny{\rm dof}}}

\newcommand{\veps}{\varepsilon}
\newcommand{\vt}{\vartheta}
\newcommand{\be}{\begin{equation}}
\newcommand{\ee}{\end{equation}}
\newcommand{\beq}{\begin{equation}}
\newcommand{\eeq}{\end{equation}}
\newcommand{\beqas}{\begin{eqnarray*}}
\newcommand{\eeqas}{\end{eqnarray*}}
\newcommand{\ep}{\varepsilon}

\newcommand{\spa}{{\mathbb R}^k}

\newcommand{\KL} {{Karh\'{u}nen-Lo\`{e}ve }}

\newcommand{\supp}{\mathrm supp}

\newcommand{\TT}{{\mathcal T}}
\newcommand{\EE}{{\mathbb E}}

\newcommand{\bbR}{\mathbb R}

\newcommand{\cG}{{\mathcal G}}
\newcommand{\cF}{{\mathcal F}}

\newcommand{\cE}{{\mathcal E}}
\newcommand{\bbN}{\mathbb N}

\newcommand{\Lipsh}{{Lipschitz }} 

\newcommand{\LN}{\lceil L/2\rceil}

\definecolor{darkred}{rgb}{.7,0,0}

\definecolor{darkgreen}{rgb}{0,0.7,0}

\definecolor{darkblue}{rgb}{0,0,0.7}

\newcommand{\bea}{\begin{eqnarray}}
\newcommand{\eea}{\end{eqnarray}}
\newcommand{\beas}{\begin{eqnarray*}}
\newcommand{\eeas}{\end{eqnarray*}}
\newcommand{\ds}{\displaystyle}

\def\IN{\mathbb{N}}

\def\IP{\mathbb{P}}
\def\IR{\mathbb{R}}

\def\cA{{\mathcal A}}
\def\cB{{\mathcal B}}

\def\cE{{\mathcal E}}
\def\cF{{\mathcal F}}
\def\cG{{\mathcal G}}

\def\cL{{\mathcal L}}
\def\cO{{\mathcal O}}
\def\cP{{\mathcal P}}

\begin{document}
\title{Complexity Analysis of Accelerated MCMC Methods for Bayesian Inversion
} 

\author[V.H. Hoang, Ch. Schwab and A.M. Stuart]{Viet Ha Hoang$^1$,
Christoph Schwab$^2$ and Andrew M. Stuart$^3$}

\address{
$^1$ Division of Mathematical Sciences, School of Physical and Mathematical Sciences, Nanyang Technological
University, Singapore, 637371}
\ead{vhhoang@ntu.edu.sg}

\address{
$^2$ Seminar for Applied Mathematics, ETH, 8092 Zurich, Switzerland}
\ead{christoph.schwab@sam.math.ethz.ch}

\address{
$^3$ Mathematics Institute, Warwick University, Coventry CV4 7AL, UK
}
\ead{a.m.stuart@warwick.ac.uk}



\begin{abstract}
{
The Bayesian approach to inverse
problems, in which the posterior probability
distribution on an unknown field is sampled
for the purposes of computing posterior
expectations of quantities of interest, 
is starting to become computationally feasible
for partial differential equation (PDE) inverse
problems.  Balancing the sources of error
arising from finite dimensional approximation  
of the unknown field, the PDE forward solution map 
and the sampling of the probability space 
under the posterior distribution
is essential for the design of efficient 
computational Bayesian methods for PDE inverse problems. 

We study Bayesian inversion
for a model elliptic PDE
with unknown diffusion coefficient.
We provide complexity analyses of several Markov Chain-Monte Carlo (MCMC)
methods for the efficient numerical evaluation of 
expectations under the Bayesian posterior distribution,
given data $\delta$.
Particular attention is given to bounds on the overall
work required to achieve 
a prescribed error level $\varepsilon$. 
Specifically,
we first bound the computational complexity of ``plain'' MCMC,
based on combining MCMC sampling 
with linear complexity multilevel solvers for elliptic PDE.
Our (new) work versus accuracy bounds show that
the complexity of this approach can be quite
prohibitive.  Two strategies for 
reducing the computational complexity  
are then proposed and analyzed: 
first, a sparse, parametric and deterministic 
generalized polynomial chaos (gpc)
``surrogate'' representation
of the forward response map of the PDE 
over the entire parameter space, and,
second, 
a novel Multi-Level Markov Chain Monte Carlo (MLMCMC)
strategy which utilizes sampling from a multilevel
discretization of the posterior and 
of the forward PDE.

For both of these strategies
we derive asymptotic bounds on work versus accuracy, 
and hence asymptotic bounds on the computational complexity 
of the algorithms. 
In particular we provide sufficient conditions 
on the regularity of the unknown coefficients
of the PDE, and on the approximation methods used, 
in order for the accelerations of MCMC
resulting from these strategies
to lead to complexity reductions over 
``plain'' MCMC algorithms
for Bayesian inversion of PDEs.
}
\end{abstract}

\maketitle

\section{Introduction} \label{sec:I}

Many inverse problems arising in PDEs
require determination of uncertain
parameters $u$ 
from finite dimensional observation data $\delta$.
We assume $u$ and $\delta$ to be related by
\be
\label{eq:basic}
\delta = \cG(u) + \vt\;.
\ee
Here $u$, which we assume to belong to
an (infinite-dimensional) function space, 
is an unknown or uncertain input (e.g. an uncertain coefficient)
to a differential equation 
and $\cG$ is the ``forward'' mapping taking 
one instance (one realization) of input $u$
into a finite and noisy set of observations.
We model these observations mathematically
as continuous linear functionals on the solution
of the governing partial differential equation. 
In \eref{eq:basic}, the parameter $\vt$ represents noise
arising when observing. For problems such as these
the Bayesian formulation \cite{kaipio2004,AndrewActa}
is an attractive and natural one,
because it allows for explicit incorporation of the
statistical properties of the observational noise,
because it admits the possibility of quantifying uncertainty
in the solution and because it allows for clear mathematical
modelling of the prior information required to account for
the under-determined nature of the inversion. 
However the Bayesian approach to such inverse problems requires
the probing of probability measures in infinite dimensional
spaces and is hence a substantial computational task in which
the space of the unknown parameter, the forward map and the 
probability space must all be approximated finite-dimensionally.
Choosing appropriate approximation techniques, and
balancing  the resulting errors so as to optimize the overall
computational cost per unit error, is thus a significant
problem with potential ramifications for the resolution of
many inverse problems arising in applications.
The purpose of this paper is to provide a detailed and
comprehensive theoretical study of this problem in
computational complexity. 

A key point to appreciate concerning our work herein 
is to understand it in the context of the existing statistics 
literature concerning the complexity of MCMC methods. This
statistics literature is focussed on the error stemming
from the central limit theorem estimate of the convergence of sample
path averages, and additionally on issues such as the ``burn-in'' 
time for the Markov chain to reach stationarity. 
Such analyses remain central in analyzing the computational
complexity of MCMC methods for inverse problems, but must be
understood, additionally, in the context of errors stemming from the finite
dimensional approximation of the unknown function and the
the forward model. The purpose of this paper is to provide
what we believe is the first complete, rigorous analysis of
computational complexity of MCMC methods for inverse problems
which balances the errors resulting from both the central limit
theorem and finite dimensional approximation.

The analysis is
necessarily quite complex as it involves balancing
errors from many different sources. 
Thus we address a specific inverse problem, namely
determination of the diffusion coefficient in an elliptic
PDE from measurements of the solution. 
We approximate the
probability space by means of an MCMC independence sampler, and we approximate
the parameter space by means of  Karhunen-Loeve expansion and
generalized polynomial chaos, and we approximate the forward
problem by finite element type methods. 
However the ideas are generic and, with substantial work,
will be transferable to other inverse problems, MCMC samplers and
approximation methodologies. 
In our analysis, we place uniform prior measures
on the unknown function $u$. This simplifies certain steps
in the present analysis; however, our results will also apply
to priors whose densities have bounded supports. 
Log-normal priors as used, for example, 
in \cite{scheichl2013} and the references therein,
require modifications in various parts 
of our analysis; this shall be dealt with elsewhere.
Finally, we concentrate on the independence sampler
as it is, currently,
the only function-space MCMC method
for which an analysis of convergence rates
which accounts
for sampling and discretization errors is available;
extensions to other MCMC algorithms are conceivable,
for example by continuing the analyses intiated in
\cite{hairer2013,vollmer2013}.

\subsection{Overview of Paper}
\label{ssec:Overview}
The methods we study incur two principal sources of error. 
First, the {\em sampling error} 
arising from estimating expectations conditional
on given data $\delta$ by sample averages of $M$ realizations
of the unknown $u$, drawn from the posterior measure 
$\rho^\delta$.
The error in doing so will decay as $M^{-\frac12}$ 
as the number $M$ of draws of $u$ tends to $\infty$.
Second, the {\em discretization error} 
arising from approximation
of the system response for each draw of $u$, 
i.e. the error of approximating $\cG(u)$.
For expository purposes, and to cover a wide range
of discretization techniques, 
we let $\eN$ denote the total 
number of degrees of freedom that are to be computed 
for evaluation of the Bayesian estimate.
We assume that the discretization error decays as $\eN^{-a}$, with
rate $a>0$ \footnote{logarithmic corrections
also occur, and will be made explicit in later sections}.
We also assume that the work per step of MCMC scales as $\eN^{b}$, with
$b>0$; thus the total work necessary for 
numerical realization of $M$ draws in the MCMC 
with discretized forward model scales as $M\eN^{b}$. 
If (as we show in the present paper for the independence sampler) 
the constant in the root mean square MCMC error 
bound of order ${\mathcal O} (M^{-\frac12})$ is 
independent of $\eN \to \infty$ 
then a straightforward calculation shows that the work
to obtain root mean square  
error $\veps$ will grow asymptotically 
as $\veps\rightarrow 0$, as $\veps^{-2-b/a}$.
The ratio $b/a$ is thus crucial to the overall computational
complexity of the algorithm. 

In this paper, 
we develop three ideas to speed up 
MCMC-based Bayesian inversion
in systems governed by PDEs.
The first idea, 
which underlies the preceding expository calculation
concerning complexity, is that MCMC methods
can be constructed whose convergence rate 
is indeed independent of the number of degrees of freedom $\eN$ 
used in the approximation of the forward map; 
the key idea here is to use MCMC
algorithms which are defined on 
(infinite dimensional) function spaces,
as overviewed in \cite{cotter2012mcmc}, and to use
Galerkin projections of the forward map 
into finite-dimensional subspaces
which employ Riesz bases in these function spaces.
We term this the {\em Plain MCMC Method}.
The second idea is that 
sparse, deterministic parametric representations 
of the forward map $u\mapsto \cG(u)$, 
as analyzed in \cite{BAS,ScCJG,CCDS11},
can significantly reduce $b$ 
by reducing computational complexity 
per step of the Markov chain, as the
sparse approximation of $\cG$ can be precomputed, 
prior to running the Markov chain, and 
simply needs to be evaluated at each step of the chain; 
this idea has been succesfully used in practice in
{\cite{marzouk2007stochastic,marzouk2009stochastic,marzouk2009dimensionality}}.  
We term the resulting algorithm the {\em gpc-MCMC Finite Element Method}. 
The third idea is that the representation of the
forward map can be truncated adaptively 
at different discretization levels 
of the physical system of interest. 
Then, we propose a multilevel Monte Carlo acceleration
of the MCMC method, in the spirit of the work of
Giles for Monte Carlo \cite{G08},
and prove that this allows  
further improvement of the computational complexity. 
The idea of extending Giles' work to Monte Carlo
Markov Chain, to obtain what we term {\em MLMCMC Methods},
is actively being considered by several other
research groups at present and we point in particular to
\cite{gruhlke2013}, where the ideas are developed in
the context of conditioned diffusions, and \cite{scheichl2013} where
uncertainty quantification in subsurface flow is studied.

To give the reader some ideas of what will come ahead, 
the following table summarizes the complexity required 
to obtain the approximation 
error $\ep$ in mean square with respect to the probability space generated by the Markov 
chains for the three methods studied
\footnote{ignoring multiplicative logarithmic factors for clarity of exposition}
\vskip 20pt
\centerline{\begin{tabular}{|c|c|c|c|}
\hline
 & Plain MCMC & gpc-MCMC & MLMCMC\\
\hline
Assume:& Assumptions \ref{assump:sumpsi}(i),(ii), \ref{assump:FEcomplexity} & 
Assumption \ref{ass:gpcplxity} & 
Assumptions \ref{assump:sumpsi}(i),(ii), \ref{assump:FEcomplexity}
\\
\hline
$\mbox{Number of} \atop\mbox{degrees of freedom}$ & $O(\ep^{-d-2})$ &  $O(\ep^{-1/\tau})$ &$O(\ep^{-d})$
\\
\hline
$\mbox{ Number of floating}\atop\mbox{ point operations}$ 
& $O(\ep^{-d-2-1/q})$ 
& $O(\ep^{-\max(\alpha/\tau,2+1/\sigma)})$ 
& $O(\ep^{-d-1/q})$ 
\\
\hline 
\end{tabular}
}
\vskip 10pt
\noindent 
Here,  $d\geq 2$ denotes the dimension of the spatial variable
in the elliptic inverse problem and $q$ determines the rate of decay in the sequence of 
coefficients of the \KL expansion of the unknown coefficient
 in Assumption \ref{assump:sumpsi}(ii).
The parameters $\tau$, $\alpha$ and $\sigma$ concern the approximation properties
of the gpc solver and are detailed in Assumption \ref{ass:gpcplxity}; $\tau$ 
quantifies the rate of convergence and $\alpha$ the polynomial scaling 
of the floating point complexity with respect to the total number of degrees of 
freedom that we assume for the gpc solver, whilst $\sigma$ relates the number of
degrees of freedom to the number of terms in the gpc approximation. 
The existence of methods which realize the assumptions made
is discussed in detail in Appendices D and E; in particular
gpc-based adaptive forward solvers of infinite-dimensional,
parametric PDEs whose overall complexity scales polynomially
in the number of degrees of freedom in the discretization of 
PDE of interest have recently become available; we refer to \cite{ScCJG,cjg:fdisc,BAS} and 
the references there. The three main results of the paper, substantiating
the displayed table, are then Theorems \ref{t:fpe}, \ref{thm:improvedMCMC}
and \ref{thm:MLMCMCerr}. 
{
Our results show that the MLMCMC FEM
can achieve a given accuracy in work 
equal asymptotically, as the accuracy $\varepsilon \to 0$,
to that of performing a single step of ``plain'' MCMC FEM;
}
they also show that the
gpc-accelerated MCMC FEM is superior 
to the ``plain'' MCMC FEM
when the compression afforded by sparse gpc FE approximations
of the parametric forward solution operator is large,
ie. 
when the parameter $\tau$ is close to $1/d$ and 
when $\alpha$ is close to 1 -- parameter regimes
which do hold in many cases 
(see, e.g, \cite{ScCJG} and the references there for 
 stochastic Galerkin discretizations of the parametric forward problem). 

The paper is organized as follows. 
Section \ref{sec:EllInvPb} is devoted to the specific
elliptic inverse problem which we study for illustrative
purposes, and includes Bayesian formulation, a discussion of
approximation techniques for the unknown field and the forward
map, and properties of the independence MCMC sampler that
we study. In section \ref{sec:S} we analyze 
the {\em Plain MCMC Method}, 
in Section \ref{sec:I2} the {\em gpc-MCMC Finite Element Method} 
and in Section \ref{sec:M} 
we study the {\em MLMCMC Method}; these sections contain
 Theorems \ref{t:fpe}, \ref{thm:improvedMCMC}
and \ref{thm:MLMCMCerr} respectively.
Section \ref{sec:Concl} contains our concluding remarks.
The paper also includes five appendices which are devoted
to mathematical developments which are necessary for the
theoretical advances made in this paper, but which employ
ideas already in the literature. In particular Appendix A concerns
Lipschitz dependence of the forward problem on the unknown function $u$, 
whilst Appendix B
is devoted to the formulation of Bayesian inverse problems
on measure spaces and Appendix C to the properties of independence
samplers in this setting. Appendices D and E concern the
finite element and polynomial chaos approximation properties,
and demonstrate how the assumptions made concerning them may
be verified.
\subsection{Overview of Notation}
\label{ssec:OvNotat}
Because of the many different approximations used in this
paper, together with the mixture of probability and PDE,
the notation can be quite complicated. For this reason
we devote this subsection to an overview of the notation, 
in order to help the reader.

To be concrete we assume that the observational noise appearing
in (\ref{eq:basic}) is a single realization of a centred Gaussian $N(0,\Sigma)$. 
If a prior measure $\rho$ is placed on the unknown $u$ then
a Bayesian formulation of the inverse problem 
leads to the problem of probing the {posterior} probability measure
$\rho^\delta$ given by
\be
\frac{d\rho^\delta}{d\rho}(u) \propto \exp(-\Phi(u;\delta))\;,
\label{eq:posterior}
\ee
where 
\be
\label{eq:phi}
\Phi(u;\delta) = {1\over 2}| \delta - \cG(u)|^2_\Sigma
\;.
\ee
Here $|\cdot|_{\Sigma}=|\Sigma^{-\frac12}\cdot|$ with
$|\cdot|$ the Euclidean norm. 
In Section \ref{sec:EllInvPb} and Appendix B
we verify that the formulae \eref{eq:posterior} and \eref{eq:phi},
which represent Bayes' formula in the infinite dimensional setting, 
hold for the model elliptic inverse problem
which is under consideration here.

Our approximation of the unknown field $u$ will be performed
through truncation of the coefficients and
Finite Element Galerkin discretizations will be used to approximate
the forward solution map $\cG$.  

On the probability side we use the following notation: 
$\cB^k$ denotes the sigma algebra of Borel subsets of $\mathbb{R}^k$.
For a probability space $(\Omega,\cA,\rho)$
consisting of the set $\Omega$ of elementary 
events, a sigma algebra $\cA$ and a probability
measure $\rho$, and a separable Hilbert space 
$H$ with norm $\|\cdot \|_H$ and for a
summability exponent $0<p\leq \infty$ 
we denote by $L^p(\Omega,\rho;H)$ the Bochner
space of strongly measurable mappings 
from $\Omega$ to $H$ which are $p$-summable.

We denote by $\IE^\mu$ the
expectation with respect to a probability measure $\mu$ on the 
space $U$ containing the unknown function $u$. 
In the following we will finite-dimensionalize both the
space $U$, in which the unknown function $u$ lies,
and the space containing the response of the forward
model. The parameter $J$ denotes the truncation level 
of the coefficient expansion \eref{eq:a} used
for the unknown function, and the parameter $l$
denotes the spatial finite element discretization level 
introduced in section \ref{sec:S}.
The parameters $N$ and $\cL$
denote the cardinality of the set 
of the chosen active gpc coefficients
and the set of finite element discretization
levels for these coefficients 
which is introduced in section \ref{sec:I2}. 
We employ multilevel FE discretizations on mesh levels $l = 0,1,...,L$ 
of the forward problem
together with multilevel approximations of the Bayesian posterior
which are indexed by $l'=0,1,\ldots,L$,
and combine these judiciously with {\em a 
discretization level dependent
sample size $M_{ll'}$} within the MCMC method.
The measure $\mu$ will variously be chosen as the
prior $\rho$, the posterior $\rho^{\delta}$, and
various approximations of the posterior
such as $\rho^{J,l,\delta}.$

We denote by 
$\cP_{u^{(0)}}$, $\cP_{u^{(0)}}^{J,l}$ and $\cP_{u^{(0)}}^{N,\cL}$ 
probability measures  
on the probability space generated 
by the MCMC processes 
detailed in the following, when 
conditioned on the initial point $u^{(0)}$ with
the acceptance probability for the
Metropolis-Hastings Markov chain being
$\alpha$ defined in \eref{eq:alpha}, 
 $\alpha^{J,l}$ defined in \eref{eq:alphaJL}, 
and 
$\alpha^{N,\cL}$ defined in \eref{eq:alphaNcL} respectively
for the problems on the
full, infinite dimensional space and its truncations. 
We then denote by $\cE_{u^{(0)}}$, $\cE_{u^{(0)}}^{J,l}$ and $\cE_{u^{(0)}}^{N,\cL}$ expectation 
with respect to $\cP_{u^{(0)}}$, $\cP_{u^{(0)}}^{J,l}$ 
and $\cP_{u^{(0)}}^{N,\cL}$ respectively.

If the initial point $u^{(0)}$ of these Markov chains is 
distributed with respect to an initial probability measure 
$\mu$ on $U$, then we denote the probability measure on the 
{probability} space that describes these Markov chains by
$\cP^\mu$, $\cP^{\mu,J,l}$ and $\cP^{\mu,N,\cL}$, 
and the corresponding expectation accordingly 
by $\cE^\mu$, $\cE^{\mu,J,l}$ and $\cE^{\mu,N,\cL}$. 
As these notations are somewhat involved, 
     for the convenience of the reader we collect them 
      in the following table.
\vskip 10pt
\begin{tabular}{|l|c||c|c|c|}
\hline
\multicolumn{2}{|c||}
{Acceptance probability} & $\alpha$ in \eref{eq:alpha} & $\alpha^{J,l}$ in \eref{eq:alphaJL} & $\alpha^{N,\cL}$ in \eref{eq:alphaNcL}\\
\hline \hline
\rule{0pt}{4ex}starting with $u^{(0)}$ & $\mbox{Probability of the}\atop\mbox{ Markov chain\ \ }$ &$\cP_{u^{(0)}}$ &$\cP_{u^{(0)}}^{J,l}$ &$\cP_{u^{(0)}}^{N,\cL}$ \\
\cline{2-5}
\rule{0pt}{4ex}& Expectation & $\cE_{u^{(0)}}$&$\cE_{u^{(0)}}^{J,l}$& $\cE_{u^{(0)}}^{N,\cL}$\\
\hline \hline
\rule{0pt}{4ex} $\mbox{starting with}\atop\mbox{distribution $\mu$}$& $\mbox{Probability of the}\atop\mbox{ Markov chain\ \ }$ & $\cP^\mu$ & $\cP^{\mu,J,l}$& $\cP^{\mu,N,\cL}$\\
\cline{2-5}
\rule{0pt}{4ex}& Expectation & $\cE^\mu$ & $\cE^{\mu,J,l}$ & $\cE^{\mu,N,\cL}$\\
\hline
\end{tabular}
\vskip 10pt

Finally, in Section \ref{sec:M}, we will work with 
the probability measure ${\mathbf P}_L$, 
on the {probability space corresponding to} 
a sequence of {independent} Markov chains created 
by the multilevel-MCMC procedure, and 
with ${\mathbf E}_L$, the expectation with respect 
to this probability measure. 
The definition of these measures will be given at the
beginning of Section \ref{sec:M}.

\section{Elliptic Inverse Problem and Approximations}
\label{sec:EllInvPb}
In this section,
we formulate the elliptic problem of
interest in Section \ref{ssec:E}, 
formulate the resulting Bayesian inverse problem
in Section \ref{ssec:BayEllInvPb}, describe the independence
sampler used to probe the posterior distribution in 
Section \ref{ssec:ind}
and convergence rate estimates of the
finite element and gpc methods in Sections \ref{ssec:FeApprFwdPbm}
and \ref{sec:SPTApprox}, respectively. 
The results in this section all follow from 
direct application of well-known theories of inverse
problems, MCMC and discretization techniques for PDEs; however, 
to render the exposition self-contained,
proofs and references are provided in Appendices A--E.
\subsection{Forward Problem}
\label{ssec:E}
Let $D$ be a {bounded \Lipsh} domain in $\IR^d$. 
For $f\in L^2(D)$, we consider the elliptic problem
\be
-\nabla\cdot(K(x)\nabla P(x))=f(x)\mbox{ in }D,\ \ 
P=0\mbox{ on }\partial D 
\;.
\label{eq:prob}
\ee
Throughout we assume that the domain $D$ is a convex polyhedron with plane sides. 
The coefficient $K \in L^\infty(D)$ in \eref{eq:prob}
is parametrized by the series expansion
\be
K(x)=\bar K(x)+\sum_{j\geq 1} u_j\psi_j(x),\quad x\in D
\label{eq:a}
\ee
where the $u_j$ are normalized to satisfy $\sup_j {|u_j|} \le 1$
and where the summation may be either infinite or finite.
In the next subsection we formulate the problem of
determining the (uncertain) function $K$, 
or equivalently the sequence $\{u_j\}_{j \ge 1}$, from a finite
noisy set of observations comprised of linear functionals
of the solution $P.$
Where it is notationally helpful to do so, we will write
$K(x,u)$ and $P(x,u)$
for the coefficient and solution of \eref{eq:prob} respectively. 

The following 
{sparsity assumptions on $K$ in \eref{eq:a},
which we will use in various combinations throughout the paper, 
imply the bounded invertibility of the parametric forward map
$U \ni u \mapsto P(\cdot,u)\in V$.
They also imply sparsity of gpc representations 
of this map and allow for
 controlling various approximation errors that arise in the sequel. 
}
\begin{assumption}\label{assump:sumpsi}
The functions $\bar{K}$ and $\psi_j$ in \eref{eq:a} 
are in $L^\infty(D)$ and:
\begin{enumerate}

\item there exists a positive constant $\kappa$ 
such that
\[
\sum_{j\geq 1} 
\|\psi_j\|_{L^\infty(D)}
\le 
{\kappa\over 1+\kappa}\bar K_{\min}, 
\]
where $\bar K_{\min}={\rm essinf}_x\bar K(x)>0$;

\item the functions $\bar K$ and $\psi_j$ in \eref{eq:a} 
are in $W^{1,\infty}(D)$ and 
there exist positive constants $C$ and $q$
such that for all $J\in \IN$,
the sequence $\{\psi_j\}$ in 
the parametric representation \eref{eq:a} of the uncertain coefficient
in \eref{eq:prob} satisfies
%
$$
\sum_{j\geq 1} \|\psi_j\|_{W^{1,\infty}(D)}
< \infty,
\quad 
\sum_{j>J}\|\psi_j\|_{L^\infty(D)}< CJ^{-q};
$$
\item  
the coefficients $\psi_j$ are arranged in
decreasing order of magnitude of $\|\psi_j\|_{L^\infty(D)}$ 
and there are constants $s>1$ and $C>0$ such that 
\[
\forall j\in \IN:\quad \|\psi_j\|_{L^\infty(D)}\le C j^{-s};
\]

\item for all $j\in \IN$, $\psi_j \in W^{1,\infty}(D)$
and there exists a constant $C>0$ such that
\[
\forall j\in \IN:\quad
\|\nabla\psi_j\|_{{L^{\infty}(D)}} \le {C j^{-s'}}
\quad\mbox{for some} \quad 1< s'\le s
\;.
\]
\end{enumerate}
\end{assumption}
{
Sparsity of the unknown diffusion coefficient $K$, 
which is expressed in terms of 
the decay rate for the coefficients $\psi_j$ of the expansion \eref{eq:a} 
}
in Assumption \ref{assump:sumpsi} (iii) and (iv) 
holds when the covariance of the random coefficient $K(x,\omega)$ 
is sufficiently smooth, as shown in Bieri et al \cite{BAS}. 
With the decay rate in Assumption \ref{assump:sumpsi}(iii), 
the constant $q$ in (ii) {can be chosen as $q = s-1$.}

We denote by $U=[-1,1]^{\IN}$ the set of all 
sequences $u=(u_j)_{j\geq 1}$ of coordinates $u_j$ 
taking values in $[-1,1]$ and note that this is 
the unit ball in $\ell^{\infty}(\IN)$.
We equip the parameter domain $U$ with the product sigma
algebra $\Theta=\bigotimes_{j=1}^\infty{\mathcal B}([-1,1])$.
Due to Assumption \ref{assump:sumpsi}(i), 
for any $u\in U$ the series \eref{eq:a}
converges in $L^\infty(D)$. 
Furthermore, it also follows from this assumption that
there exist finite positive constants 
$K_{\min}$ and $K_{\max}$ such that, for all $(x,u)\in D \times U$,
\beq\label{eq:KminKmax}
K_{\min}:=\bar K_{\min}/(1+\kappa)
\le 
K(x)
\le 
K_{\max}:={\rm esssup}_x\bar K(x)+\kappa\bar K_{\min}/(1+\kappa).
\eeq

We let $V = H^1_0(D)$, whilst $V^*$ denotes its dual space. 
{
We equip $V$ with the norm 
$\| P \|_V = \| \nabla P \|_{L^2(D)}$.
By \eref{eq:a}, 
$K(x)$ is bounded below uniformly with respect
to $(x,u) \in D \times U$ and, 
for every $u\in U$,
\[
\begin{array}{rcl}
\fl\ds
K_{min} 
\| P(\cdot,u)\|^2_V
&=& \ds 
K_{min} 
(\nabla P(\cdot,u),\nabla P(\cdot,u))
\leq 
(K(\cdot,u) \nabla P(\cdot,u),\nabla P(\cdot,u)) 
\\
&=& \ds 
(f, P(\cdot,u))
\leq { \| f \|_{V^*} }\| P(\cdot,u) \|_V 
\;,
\end{array}
\]
where $(\cdot,\cdot)$ denotes the inner product in $L^2(D)$ and $(L^2(D))^d$.
It follows that 
\be\label{eq:stndLaxMilg}
\sup_{u\in U} \| P(\cdot,u) \|_V 
\leq 
\frac{ \| f \|_{V^*} }{ K_{min} } 
\;.
\ee
}
In fact we have the following, proved in Appendix A.
\begin{proposition}\label{prop:measurable}
Under Assumption \ref{assump:sumpsi}(i), the solution 
$P: U \mapsto V = H^1_0(D)$ is \Lipsh when viewed
as a mapping from the unit ball in $\ell^{\infty}(\IN)$
to $V$.
It is in particular measurable,
as a mapping from the measurable space 
$(U,\Theta)$ to $(V,\cB(V))$.
\end{proposition}
%
%
\hfill$\Box$
\subsection{Bayesian Elliptic Inverse Problem}
\label{ssec:BayEllInvPb}
We now define the Bayesian inverse problem. 
For $\cO_i\in V^*$, $i=1,\ldots,k$, 
which denote $k$ continuous, linear ``observation'' functionals on $V$, 
we define a map $\cG:U\to {\mathbb R}^k$ as 
\be\label{eq:defG}
U\ni u\mapsto \cG(u) 
:= 
(\cO_1(P(\cdot,u)),\cO_2(P(\cdot,u)),\ldots,\cO_{k}(P(\cdot,u)))
\in \bbR^{k}
\;.
\ee
In \eref{eq:basic},
by $\vartheta$ we denote an observational noise 
which is statistically independent of the input $u$ and which
we assume to be Gaussian, ie. with
distribution $N(0,\Sigma)$ in ${\mathbb R}^k$, 
with positive definite covariance matrix $\Sigma$.
We model the noisy 
observed data $\delta$ for $\cG(u)$ by
\[
\delta = \cG(u)+\vartheta \;.
\]
For Bayesian inversion, we place a prior 
measure on $u$ by assuming that $u_j:\Omega\to [-1,1]$ 
comprises a sequence of 
independent random variables $u_j:\Omega\to [-1,1]$
in the series expansion \eref{eq:a}.
On the measurable space $(U,\Theta)$ defined above 
we define the countable product 
probability measure
\[
\rho=\bigotimes_{j\geq 1}{du_j\over 2},
\]
where $du_j$ is the Lebesgue measure on $[-1,1]$. 
As $u_j$ are uniformly distributed on $[-1,1]$, 
the measure $\rho$ is the law of 
the random vector $u=(u_1,u_2,\ldots)$ in $U$. 
As the random variables $u_j(\omega)$ were 
assumed independent,
the probability measure on realizations of 
random vectors $u\in U$ is a product measure:
for $S=\prod_{j\geq 1} S_j$,
\[
\rho(S)=\prod_{j\geq 1} \IP(\{\omega:u_j\in S_j\})
\;.
\]
Combining the prior and likelihood, the posterior measure on $u$
given $\delta$ can be explicitly  determined, using the general
theory from Appendix B, and a stability/well-posedness
estimate;
the proof of the following result 
{for the model problem} 
is provided there.
\begin{proposition} \label{p:wp}
Let Assumption \ref{assump:sumpsi}(i) hold.
The conditional probability measure 
$\rho^\delta(du) = \IP(du|\delta)$ on $U$ 
is absolutely continuous with respect to $\rho(du)$ 
and satisfies
\[
{d\rho^\delta \over d\rho}\propto \exp(-\Phi(u;\delta))\;.
\]
Furthermore, 
for every $r>0$ for every $\delta,\delta'$ such that 
$|\delta|_\Sigma ,|\delta'|_\Sigma \le r$,
there exists $C=C(r)>0$ such that
$$ 
d_{\rm{Hell}}(\rho^\delta,\rho^{\delta'}) \le C(r)|\delta - \delta'|_\Sigma
\;
$$
where $d_{\rm{Hell}}$ denotes the Hellinger distance of the 
measures $\rho^\delta,\rho^{\delta'}$ (see, eg., 
equation \eref{eq:Hell} in Appendix B).
\end{proposition}

\begin{remark} 
{The proof}
of the preceding proposition
shows that $|\cG(u)|$ is uniformly bounded for $u$
in $U$. As a consequence there exists $\eM(r)>0$ which is a uniform
bound on $\Phi(u;\delta)$ for all $|\delta| \le r$ 
and for all $u \in U$.
This bound is also uniform with respect to truncation of
the infinite series \eref{eq:a} for $K$, since
this corresponds to a particular choice of some of 
the coefficients of $u \in U$, and with
respect to finite element approximation of the solution
of \eref{eq:prob}, since the uniform upper bound
on $|\cG(u)|$ {will be preserved {under Galerkin projections
of the elliptic} problem \eref{eq:prob} into finite-dimensional subspaces
$V_h \subset V$ (or, more generally, under any stable discretization
of the forward problem of interest).
}
\label{r:uniform2}
\end{remark}

\subsection{Independence Sampler}
\label{ssec:ind}

To approximate expectations with respect to the posterior
measure $\rho^{\delta}$ constructed in the previous section
we will use MCMC methods and, in particular, the independence
sampler.  To this end we define, for any $u,v\in U$,
\begin{equation}
\label{eq:alpha}
\alpha(u,v)=1\wedge \exp\bigl(\Phi(u,\delta) - \Phi(v,\delta)\bigr)
\;.
\end{equation}
The Markov chain 
$\{u^{(k)}\}_{k=1}^{\infty} \subset U$ 
is {then} constructed as follows:
given the current state $u^{(k)}$, 
we draw a proposal $v^{(k)}$
independently of $u^{(k)}$ from the
prior measure $\rho$ appearing in \eref{eq:posterior}. 
Let $\{w^{(k)}\}_{k \ge 1}$ denote an i.i.d sequence 
with $w^{(1)} \sim {\mathcal U}[0,1]$ and with
$w^{(k)}$ independent of both $u^{(k)}$ and $v^{(k)}$.
We then determine the next state $u^{(k+1)}$ via the formula 
\be
\label{eq:MC}
\fl u^{(k+1)}
=
{\bf 1}\bigl(\alpha(u^{(k)},v^{(k)})\ge w^{(k)}\bigr)
v^{(k)}
+\Bigl(1-{\bf 1}\bigl(\alpha(u^{(k)},v^{(k)})\ge w^{(k)}\bigr)
\Bigr)u^{(k)}
\;.
\ee
Thus we choose to move from $u^{(k)}$ to $v^{(k)}$ with probability
$\alpha(u^{(k)},v^{(k)})$, and to remain at $u^{(k)}$ with 
{probability $1-\alpha(u^{(k)},v^{(k)})$. The following is proved
in Appendix C.

\begin{proposition} \label{t:mc}
Let Assumption \ref{assump:sumpsi}(i) hold. 
Equation \eref{eq:MC} 
defines a Markov chain $\{u^{(k)}\}_{k=0}^{\infty}$ which is reversible
with respect to $\rho^\delta$. If
$p(u,A)$ denotes the transition kernel for the
Markov chain, and $p^n(u,A)$ its $n^{th}$ iterate, 
then there is $R \in (0,1)$ such that, for all $n\in \bbN$,
$$
\|p^n(u,A)-\rho^\delta\|_{{\rm {\tiny TV}}} 
\le 
2(1-R)^n\;.
$$
For 
every bounded, continuous 
function $g:U \to \bbR$,
there holds, $\cP_{u^{(0)}}$ almost surely,
\begin{equation}
\label{eq:ex}
\frac{1}{M}\sum_{k=1}^M g(u^{(k)})=\IE^{\rho^\delta}[g(u)] + c\xi_M M^{-\frac12}
\end{equation}
where $\xi_M$ is a sequence of random variables
which converges weakly as $M \to \infty$ to 
$\xi \sim N(0,1)$ and where $c$ is a positive constant
which is independent of $M$ (it depends only on $\eM(r)$ in 
Remark \ref{r:uniform2} and on $\sup_{u \in U}|g(u)|$).
Furthermore, we have the mean square error bound:
there exists $C>0$ such that for every bounded $g:U\mapsto \IR$ and
every $M\in \IN$
\[
\Big(\cE^{\rho}\Big[\Big|\IE^{\rho^\delta}[g(u)]
-
\frac{1}{M}\sum_{k=1}^M g(u^{(k)})\Big|^2\Big]\Big)^{1/2}
\le 
C{\sup_{u\in U}|g(u)|} M^{-1/2}\;.
\]
\end{proposition}

\begin{remark}
\label{r:uniform}
The proof of this proposition uses fairly standard methods
from the theory of Markov chains. 
In our context a key observation regarding the proof is that all
constants  -- in particular the constant $R$ -- 
depend only on the upper bound on the negative log likelihood
$\Phi$, and on the supremum of $g$. 
Hence, if we can show for {finite element and \KL}
approximations of the forward map $\cG$
in the physical domain $D$ that
these approximations are such that the 
potential $\Phi$ is stable under these approximations,
then the conclusions of the
preceding theorem will hold with constants 
that are uniformly bounded 
with respect to all approximation parameters.

Note also that our choice of uniform priors and the independence
sampler means that the issue of``burn-in'' does not play a significant
role in our analysis; in particular the total variation metric
convergence bound is independent of initialization. When
generalizing our work to other priors and other MCMC methods ``burn-in''
effects may become more pronounced in the analysis.
\end{remark}

\subsection{Finite Element Approximation of the Forward Problem}
\label{ssec:FeApprFwdPbm}
Assumptions \ref{assump:sumpsi}(i),(ii) are imposed
throughout what follows regarding the finite element method. Thus, from \eref{eq:a}, 
we deduce that $K(\cdot,u)\in W^{1,\infty}(D)$ for all $u\in U$.  
We describe an approximation of the forward problem
based on finite element representation of the solution
$P$ of \eref{eq:prob}, together with truncation
of the series \eref{eq:a}. We start by discussing
the finite element approximation. 
Recalling that the domain $D$ is a bounded \Lipsh polyhedron 
with plane sides, we denote by $\{{\mathcal T}^l\}_{l=1}^\infty$ 
the nested sequence of simplices which is defined inductively as follows: 
first we {subdivide} $D$ 
into a regular family ${\mathcal T}^0$ of simplices $T$; 
for $l\geq 1$, the regular simplicial mesh ${\mathcal T}^l$ in $D$
is obtained by recursive subdivision, i.e. 
each simplex in ${\mathcal T}^{l-1}$ is divided into 
$2^d$ subsimplices.
Then, the meshwidth $h_l := \max\{{\rm diam}(T): T\in \TT^l \}$ 
of  ${\mathcal T}^l$ is $h_l = 2^{-l}h_0$.
Based on these triangulations,
we define a nested multilevel family of spaces of continuous, 
piecewise linear functions on ${\mathcal T}^l$ as
\[
V^l=\{u\in V:\ u|_T\in{\mathcal P}_1(T)\ \forall\,T\in {\mathcal T}^l\},
\]
where ${\mathcal P}_1(T)$ denotes the set of linear polynomials in 
the simplex $T\in {\mathcal T}^l$. 
Approximating the solution of the parametric, deterministic
problem \eref{eq:prob} from the finite element spaces $V^l$
introduces a {\em discretization error} which is well-known 
to be bounded by the {\em approximation property} of the $V^l$:
there exists a positive constant $C>0$ which is independent of $l$
such  that for all $P\in H^2(D) \cap H^1_0(D)$ 
and for every $0 < h_l\leq 1$ holds
\beq\label{eq:FEApproxrate}
\inf_{Q\in V^l}\|P-Q\|_{V} \le C h_l\|P\|_{H^2(D)},
\eeq
where $h_l = O(2^{-l})$
and where the constant $C>0$ depends only on 
$\TT^0$.

We now discuss the effect of dimensional truncation,
ie. of truncating the infinite series 
\eref{eq:a} for the unknown diffusion coefficient $K$ 
in problem \eref{eq:prob} after $J$ terms, 
as
\beq\label{eq:KtruncJ}
K^J(x,u)=\bar K(x)+\sum_{j=1}^J u_j\psi_j(x)\, x\in D,\, u\in U
\;.
\eeq
We thus consider the parametric, 
deterministic family of approximate {elliptic} problems
\be
-\nabla\cdot(K^J(\cdot,u)\nabla P^J(\cdot,u)) = f,\ \ 
P^J=0\mbox{ on }\partial D 
\;.
\label{eq:probJ}
\ee
Under Assumptions \ref{assump:sumpsi}(i) and (ii),
\eref{eq:Lipschitzofp1} shows that there exists a constant $C>0$ such that for all
$J\in \IN$ and all $u\in U$
\begin{equation}
\begin{array}{rcl}
\fl 
\|P(\cdot,u)-P^J(\cdot,u)\|_{V} 
&\le& \ds 
C\|P(\cdot,u)\|_{V}\|K(\cdot,u)-K^J(\cdot,u)\|_{L^\infty(D)}
\\
&\le& 
C\|P(\cdot,u)\|_{V} J^{-q} 
\leq 
\frac{C}{K_{min}} J^{-q} \| f \|_{V^*} 
\;.
\end{array}
\label{eq:probJ2}
\end{equation}

We consider the finite element
approximation of the truncated problem \eref{eq:probJ}: 
given $J,l\in \IN$, 
find $P^{J,l}(\cdot,u)\in V^l$ such that
for all $\phi \in V^l$
\be
\int_DK^J(x,u)\nabla P^{J,l}(x,u)\cdot\nabla\phi(x)dx
=
\int_D f(x)\phi(x)dx
\;.
\label{eq:fe}
\ee
We make the following assumption on the complexity of solving this discrete equation
the justification of which is provided in Appendix D.
\begin{assumption}\label{assump:FEcomplexity}
{For $J\in \mathbb{N}$ as in \eref{eq:KtruncJ},}
the solution $P^J(\cdot,u)$ of \eref{eq:probJ} is uniformly bounded in 
$W:=H^2(D)\cap H^1_0(D)$ {with respect to $J\in \mathbb{N}$ and to $u\in U$}. 
The {matrix of the Galerkin approximated problem \eref{eq:fe}}
has $O(l^{d-1}2^{dl})$ non-zero entries and has a uniformly 
bounded condition number for all $J, l$ and $u$. 
{There exists $C>0$ such that, for all $J\in \mathbb{N}$ 
and for all $u\in U$, 
the FE error for the approximating problem \eref{eq:fe} satisfies
}
\begin{equation}
\|P^{J}(\cdot,u)-P^{J,l}(\cdot,u)\|_{V}
\le 
C 2^{-l}\|P^J(\cdot,u)\|_{W}\;.
\label{eq:fe2}
\end{equation}
\end{assumption}

%
From this assumption, we obtain the following error bound.
\begin{proposition}
\label{p:totapp}
Consider {the} approximation of the elliptic problem
\eref{eq:prob} via 
finite element solution of the truncated problem \eref{eq:probJ}, 
under
Assumptions \ref{assump:sumpsi} (i), (ii)
and Assumption \ref{assump:FEcomplexity}. 
Then there exists a constant $C>0$ such that 
for every $J,l\in \bbN$ and for every $u\in U$
it holds that
the Finite Element solutions
$P^{J,l}(\cdot,u)$ are uniformly $V$-stable
in the sense that, for every $J,l\in \mathbb{N}$
and for every $u\in U$ there holds
\be\label{eq:PJlstab}
\sup_{J,l\in \IN}\sup_{u\in U}
\| P^{J,l}(\cdot,u) \|_V \leq \frac{C}{K_{min}} \| f \|_{V^*}
\;.
\ee
Moreover, there exists $C>0$ such that for every $u\in U$
holds the error bound
\be\label{eq:PJlerr}
\|P(\cdot,u)-P^{J,l}(\cdot,u)\|_{V} 
\le 
C(2^{-l}\|P^J(\cdot,u)\|_{W} + J^{-q}\|P(\cdot,u)\|_{V})
\;.
\ee
\end{proposition}
\subsection{Sparse Tensor gpc-Finite Element {Surrogate}
            of the Parametric Forward Problem} 
\label{sec:SPTApprox}
By \eref{eq:stndLaxMilg},
the solution $P(\cdot,u)$ of problem \eref{eq:prob} 
is uniformly bounded in $V$ by $\|f\|_{V^*}/K_{\min}$. 
Therefore, from Proposition \ref{prop:measurable}, 
we deduce that $P(\cdot,\cdot)\in L^2(U,\rho; V)$. 
Thus
the parametric solution admits a polynomial
chaos type representation in $L^2(U,\rho; V)$.
To define it, we 
denote by $L_n(u_n)$ the Legendre polynomial of degree $n$,
normalized such that
\[
{1\over 2}\int_{-1}^1 | L_n(\xi)|^2 d\xi = 1\;.
\]
By ${\mathcal F}$ we denote the set of all sequences 
$\nu=(\nu_1,\nu_2,\ldots)$ of nonnegative integers
such that {only finitely many $\nu_j$ are nonzero.} 
We define
\begin{equation}\label{eq:Lnu}
L_\nu(u) = \prod_{j\geq 1} L_{\nu_j}(u_j)\;.
\end{equation}
Since $L_0 \equiv 1$, for each $\nu\in \cF$
the product contains only finitely many nontrivial factors. 
The set $\{ L_\nu : \nu\in{\mathcal F} \}$ 
forms an orthonormal basis for $L^2(U,\rho)$. 
We can therefore expand $P(\cdot,u)$ into the 
Legendre expansion
\[
P(\cdot,u) = \sum_{\nu\in {\mathcal F}}P_\nu(\cdot)L_\nu(u),
\]
where $P_\nu := \int_U P(\cdot,u)L_\nu(u) d\rho(u) \in V$. 
By the $L^2(U,\rho)$ orthonormality of the set $\{ L_\nu: \nu\in \cF\}$, 
Parseval's equation in the Bochner space $L^2(U,\rho;V)$ 
takes the form
$$
\forall P\in L^2(U,\rho;V):
\quad 
\| P \|_{L^2(U,\rho;V)}^2 
= 
\sum_{\nu\in \cF} \| P_\nu \|_V^2 \;.
$$
{
We now define the sparse tensor surrogate forward maps 
where we only use a finite (sparse) subset of the gpc mode index set $\cF$.
For any integer $N$, choose an index set 
$\Lambda_N\subset \cF$ of gpc coefficients $P_\nu\in V$
which are to be included into the surrogate map, 
subject to the constraint $\#(\Lambda_N) \leq N$, 
and a second set 
${\cal L}(\Lambda_N) := (l_\nu)_{\nu\in\Lambda_N}\subset \IN^N$ 
of Finite-Element discretization levels
for the FE approximation of the active $P_\nu$, $\nu\in \Lambda_N$.
Then we consider surrogate forward maps
$P_{N,\cL}$ for the response $P$ 
which are of the form
\be
P_{N,\cL}(x,u)=\sum_{\nu\in\Lambda_N}P_{\nu,\cL}(x)L_\nu(u),\ \ P_{\nu,\cL}\in V^{l_\nu}
\;.
\label{eq:PNL}
\ee
We wish to find the sets $\Lambda_N$ and ${\cal L}(\Lambda_N)$ 
that give the best 
(or the quasi-best) approximations among all finite subsets 
$\Lambda \subset \cF$ and of ${\cal L}(\Lambda)$ subject to
a constraint on their {\em combined} cardinality
$N_{dof} = \sum_{\nu\in \Lambda_N} l_\nu$.
Convergence rates, in terms of $N_{dof}$, of 
such approximations of parametric forward maps 
have been derived recently, for example, in
\cite{CDS1}, \cite{CDS2}, \cite{HSparabolic}, \cite{HSwave}, \cite{ScCJG}, \cite{Nisto782}.
In these references, 
the following assumption has been verified for various problem classes.
}
\begin{assumption}\label{ass:gpcplxity}
There are positive constants $\sigma$, $\tau$, $\alpha$ and $\beta$ such that for each integer $N$,  
with a total budget of 
$\eN=O(N^{\sigma/\tau})$ degrees of freedom, 
a subset $\Lambda_N\subset \cF$ 
of cardinality not exceeding an absolute multiple of $N$ 
and such that, for every $\nu\in \Lambda_N$,  $|\nu|=O(\log N)$
and a {surrogate} gpc FE approximation 
$P_{N,{\mathcal L}}$ of the parametric forward solution $P$
of \eref{eq:probJ} of the form \eref{eq:PNL}
with  rate of convergence
\[
\|P-P_{N,{\mathcal L}}\|_{L^2(U,\rho;V)} 
\le 
CN_{dof}^{-\tau},
\]
can be found in 
$O(N_{dof}^\alpha(\log N_{dof})^\beta)$ 
float point operations.
\end{assumption} 
{For the model elliptic inverse problem of determining $K$
given by \eref{eq:a} from linear functionals defined
on solutions of \eref{eq:prob}}, 
subject to the sparsity conditions in
Assumption \ref{assump:sumpsi}, the preceding assumption is
verified in Appendix E; there also 
bounds for $\tau$, the rate 
$\sigma$ and for the exponents
$\alpha > 0$ and $\beta \geq 0$ in Assumption 
\ref{ass:gpcplxity} are derived, in terms of the sparsity assumptions 
on the unknown coefficient $K$ in \eref{eq:KtruncJ}.
\section{Plain MCMC}
\label{sec:S}
We study the computational complexity of the 
Plain MCMC method \eref{eq:MC} to 
generate samples from the posterior probability 
measure $\rho^\delta$ determined in the previous section.
The main results of this analysis are summarized in
Theorem \ref{t:fpe}.
The complexity analysis of the MCMC algorithm 
is of independent interest as it utilizes the emerging
idea of MCMC methods with dimension independent 
mixing properties \cite{cotter2012mcmc,hairer2013,vollmer2013}.
Furthermore, the results in the present
section will be the foundation for several accelerations
of the Plain MCMC algorithm which will be presented in Sections \ref{sec:I2} and
\ref{sec:M} ahead.  In order to obtain a constructive version of the MCMC
algorithm, we will approximate solutions of
the forward problem \eref{eq:prob} 
by applying the Finite Element Method in the physical domain $D$
to its parametric version \eref{eq:prob}
and by truncation of the expansion of
the diffusion coefficient $K$ given by \eref{eq:a}. 

Assumptions \ref{assump:sumpsi}(i),(ii) are imposed
throughout what follows.
From Assumption \ref{assump:sumpsi}(ii) and from 
\eref{eq:a}, 
we deduce that $K(\cdot,u)\in W^{1,\infty}(D)$ for all $u\in U$.  
\subsection{FE Approximation of the Posterior Measure}
We denote the vector of observables from the 
Galerkin discretized, gpc-parametrized approximate
forward solution map by
\be\label{eq:defGJl}
\cG^{J,l}(u)
=
(\cO_1(P^{J,l}(u)),\ldots,\cO_k(P^{J,l}(u))):U\mapsto \bbR^k 
\ee
and define the corresponding Bayesian potential
\begin{equation}
\label{eq:phi2}
\Phi^{J,l}(u;\delta)={1\over 2}|\delta-\cG^{J,l}(u)|^2_\Sigma \; .
\end{equation}
We define an
approximate conditional 
posterior probability measure $\rho^{J,l,\delta}$ 
on the measurable space $(U,\Theta)$ as
\[
{d\rho^{J,l,\delta}\over d\rho}\propto\exp(-\Phi^{J,l}(u;\delta))\;.
\]
Then, the 
measure $\rho^{J,l,\delta}$ is an approximation of
the Bayesian posterior  $\rho^{\delta}$
which, due to the discretization and the truncation of the expansion \eref{eq:a}
incurs an approximation error.
We now show that this error in the posterior measure 
is bounded in the Hellinger metric 
with respect to $J$ and $l$ in the same way 
as the forward error in Proposition \ref{p:totapp}. 
The proof of the following Proposition is based on a 
generalization of the method introduced in \cite{Cotteretal}.
\begin{proposition}\label{prop:approxmu}
Let Assumptions \ref{assump:sumpsi}(i),(ii) and \ref{assump:FEcomplexity} hold.
If the domain $D$ is convex and if $f\in L^2(D)$
then there exists 
a positive constant $c$ depending 
only on the data $\delta$ such that
for every $J$ and $l$ there holds
\[
d_{\rm{Hell}}(\rho^\delta,\rho^{J,l,\delta})
\le 
C(\delta)(J^{-q}+2^{-l})
\| f \|_{L^2(D) 
}
\;.
\]
\end{proposition}

{\it Proof}\ \ We denote the normalizing constants as
\[
Z(\delta) = \int_U \exp(-\Phi(u;\delta))d\rho(u),
\quad
Z^{J,l}(\delta) = \int_U\exp(-\Phi^{J,l}(u;\delta))d\rho(u)
\;.
\]
We then estimate
\beqas 
\fl 2d_{\rm{Hell}}(\rho^\delta , \rho^{J,l,\delta})^2
\\
\fl=
\int_U\Big(
Z(\delta)^{-1/2}
\exp\Big(-{1\over2}\Phi(u;\delta)\Big) 
- 
(Z^{J,l}(\delta))^{-1/2} \exp\Big(-{1\over2}\Phi^{J,l}(u;\delta)\Big)\Big)^2
d\rho(u)
\\
\le 
I_1+I_2,
\eeqas
where we defined
\beqas
I_1 
:= 
{2\over Z(\delta)} 
\int_U
\Big(\exp\Big(-{1\over2}\Phi(u;\delta)\Big)-\exp\Big(-{1\over2}\Phi^{J,l}(u;\delta)\Big)\Big)^2
d\rho(u),
\\
I_2 := 2|Z(\delta)^{-1/2} - Z^{J,l}(\delta)^{-1/2}|^2 
\int_U\exp(-\Phi^{J,l}(u;\delta))d\rho(u)
\;.
\eeqas
{We estimate $I_1$ and $I_2$. To bound $I_1$,}
given data $\delta$, for every $u\in U$ holds
\begin{eqnarray}
\Big|\exp\Big(-{1\over2}\Phi(u;\delta)\Big) 
- 
\exp\Big(-{1\over2}\Phi^{J,l}(u;\delta)\Big)\Big|
\le 
{1\over 2}|\Phi(u;\delta) - \Phi^{J,l}(u;\delta)|
\nonumber
\\
\le 
C(2|\delta| + |\cG(u)| + |\cG^{J,l}(u)|) 
|\cG(u)-\cG^{J,l}(u)|
\;.
\label{eq:LipschitzPhi}
\end{eqnarray}
Moreover, by Proposition \ref{p:totapp},
there exists a constant $C>0$ 
independent of $J$ and of $l$ such that,
for all $u\in U$, 
there holds
\beqas
|\cG(u)-\cG^{J,l}(u)|
&\le& 
C\max \{ \|\cO_i\|_{V^*}\}\|P(\cdot,u)-P^{J,l}(\cdot,u)\|_{V}
\\
&\le& \ds
C(2^{-l}\|P^J(\cdot,u)\|_{{W}} + J^{-q}\|P(\cdot,u)\|_{V})
\;.
\eeqas
By \eref{eq:stndLaxMilg} and Proposition \ref{prop:H2},
$\| P(\cdot,u)\|_V$ and $\| P^{J}(\cdot,u)\|_W$ 
are uniformly bounded with respect to $u\in U$.
{
Therefore, 
there exists a constant $C(\delta)>0$ depending
only on the data $\delta$ such that 
for all $J\in \mathbb{N}$
}
\[
\begin{array}{rcl}
I_1 &\le & \ds 
C(\delta)\EE^{\rho}(2^{-l} \| P^J(\cdot,u)\|_{W} 
+ 
J^{-q}\|P(\cdot,u)\|_{V})^2
\\
& \le & \ds
C(\delta) ( J^{-2q} \|f\|_{V^*}^2 + 2^{-2l} \| f \|_{L^2(D)}^2) 
\;.
\end{array}
\]
To estimate term $I_2$, we observe that there is a positive constant 
$C>0$ such that for every $J,l\in \IN$ holds
\[
|Z(\delta)^{-1/2}-Z^{J,l}(\delta)^{-1/2}|^2 
\le 
C(Z(\delta)^{-3}\vee Z^{J,l}(\delta)^{-3})|Z(\delta) - Z^{J,l}(\delta)|^2\;.
\]
We note that
\[
\begin{array}{rcl}
|Z(\delta) - Z^{J,l}(\delta)| 
&\le&\ds
\int_U |\exp(-\Phi(u;\delta)) - \exp(-\Phi^{J,l}(u;\delta))|d\rho(u)
\\
& \le & \ds
\int_U |\Phi(u;\delta) - \Phi^{J,l}(u;\delta)|d\rho(u)\;.
\end{array}
\]
Therefore, 
as $Z(\delta)$ and $Z^{J,l}(\delta)$ are uniformly bounded 
below for all $\delta$,
{analysis similar to that for $I_1$ shows that} 
\[
I_2 \le C(\delta)(2^{-2l}+J^{-2q}) {\| f \|_{L^2(D)}^2}.
\]
Thus
\[
d_{\rm{Hell}}(\rho^\delta , \rho^{J,l,\delta}) 
\le 
C(\delta)(2^{-l}+J^{-q}) {\| f \|_{L^2(D)}}
\;.
\]
\hfill$\Box$
\subsection{Computational Complexity of Plain MCMC}
Given $J,l \in \bbN$, and data $\delta$,
we use the MCMC method \eref{eq:MC} to sample 
the probability measure $\rho^{J,l,\delta}$. 
In so doing
we create a method for approximating integrals
of functions $g:U \to \IR$ with respect to $\rho^{\delta}$.
We use the following notation for the empirical
measure generated by the Markov chain designed to
sample $\rho^{J,l,\delta}$:
\be\label{eq:defEMrhoJldel}
E_M^{\rho^{J,l,\delta}}[g]:={1\over M}\sum_{k=1}^Mg(u^{(k)}),
\ee
where the Markov chain 
${{\mathcal C}_{J,l}
=
\{u^{(k)}\}_{k \in \IN_0}}\subset\mathbb{R}^{J}$
is started in the restriction of 
$u^{(0)}$ to $\mathbb{R}^{J}$.
It depends on the discretization level $l$ and the truncation level $J$
since it
is generated from the process \eref{eq:MC} 
with the acceptance probability \eref{eq:alpha} 
being replaced by
\be
\alpha^{J,l}(u,v)
=
1\wedge\exp(\Phi^{J,l}(u;\delta)-\Phi^{J,l}(v;\delta))
\;,\quad 
(u,v) \in U\times U
\;.
\label{eq:alphaJL}
\ee
Given $M\in\IN$ we wish to estimate the MC sampling error
\be
\IE^{\rho^\delta}[g]-E_M^{\rho^{J,l,\delta}}[g]\;.
\label{eq:ge}
\ee
We develop in the following two types 
of error bounds as $M\rightarrow \infty$ 
for \eref{eq:ge}: a 
probabilistic error bound
for {$\cP_{u^{(0)}}^{J,l}$} almost every 
realization of the Markov chain and a mean square bound.
We refer to the table in Section \ref{ssec:OvNotat} for
the notation in the next result.
\begin{proposition}\label{prop:KJh}
Let Assumptions \ref{assump:sumpsi}(i),(ii) and \ref{assump:FEcomplexity} hold.
Let $g:U \to \IR$ be a bounded  {continuous}
function on $U$ with respect to the supremum norm. 
Then, for every initial condition $u^{(0)}$ 
and for $\cP_{u^{(0)}}^{J,l}$-almost every 
realization of the Markov chain holds the error bound
\begin{equation*}
\Bigl|
\IE^{\rho^\delta} g(u)
-
E_M^{\rho^{J,l,\delta}}[g]
\Bigr| 
\le
c_1 M^{-\frac12}+c_2\bigl(J^{-q}+2^{-l}\bigr)
\end{equation*}
where $c_1 \le c_3|\xi_M|$;
$\xi_M$ is a random variable (on the probability
space generating the randomness within the Markov chain)
which converges weakly as $M \to \infty$ to
$\xi \sim N(0,1)$ and $c_2, c_3$ are 
non-random constants independent 
of $M,J$ and $l$.

Moreover, there exists a constant
$c_4$ (which is deterministic and 
depends only on the data $\delta$, 
and which is, in particular, independent of $M,J$ and $l$)
such that 
%
\beq\label{eq:MCMCmeansqconv}
\Big(\cE^{\rho,J,l}\Big[\left|
\IE^{\rho^\delta}[g] -  E_M^{\rho^{J,l,\delta}} [g]
\right|^2\Big]\Big)^{1/2}
\le
c_4( M^{-1/2}+J^{-q} + 2^{-l}).
\eeq
\end{proposition}
{\it Proof}\ \ 
As $g$ is bounded, we have from Proposition
\ref{prop:approxmu} and properties of the Hellinger
metric (specifically, from (2.7) in \cite{Cotteretal})
for every $u\in U$ that
\beq\label{eq:Hellingineq}
|\IE^{\rho^\delta} g(u)-\IE^{\rho^{J,l,\delta}}g(u)|
\le 
\bar{c}(g)
d_{\rm{Hell}}(\rho^\delta ,\rho^{J,l,\delta })
\le \bar{c}(g) C(\delta) (J^{-q}+2^{-l})
\;.
\eeq
Here, $C(\delta)$ is as in Proposition \ref{prop:approxmu}
and $\bar{c}(g)$ depends on the supremum of $g(u)$ over $U$,
but is independent of $J,l$.
By Proposition \ref{t:mc} (and Remarks \ref{r:uniform} and
\ref{r:uniform2}) we deduce the existence of a constant
$C>0$, independent of $M,J$ and $l$, such that
there holds, as $M\rightarrow \infty$, 
$\cP_{u^{(0)}}^{J,l}$ almost surely,

\[
|\IE^{\rho^{J,l,\delta }}g - {1\over M}\sum_{k=1}^M g(u^{(k)})|
\le 
C |\xi_M| M^{-1/2}
\]
where $\xi_M$ converges weakly as $M \to \infty$ to $\xi \sim N(0,1)$.
Combining this with \eref{eq:Hellingineq} gives the first assertion.

To prove the mean square error bound \eref{eq:MCMCmeansqconv}, 
we define
\[
\bar g(u) := g(u)-\IE^{\rho^{J,l,\delta}}[g],
\]
as in the proof of Proposition \ref{t:mc}.
We note that $\bar g$ depends
on $J,l$ via dependence on the approximate posterior measure $\rho^{J,l,\delta}$,
but we do not indicate this dependence explicitly. However we will use 
uniform boundedness of $\bar g$ with respect to these parameters in an essential way
in what follows.

Due to the invariance of the stationary measures $\rho^{J,l,\delta}$, 
we may write
\beqas
\fl{1\over M}\cE^{\rho^{J,l,\delta},J,l}\Big[\Big|\sum_{k=1}^M\bar g(u^{(k)})\Big|^2\Big]
&=&\ds
\IE^{\rho^{J,l,\delta}}[\bar g(u^{(0)})^2]
+
2{1\over M}\sum_{k=1}^M\sum_{j=k+1}^M\cE^{\rho^{J,l,\delta},J,l}[\bar g(u^{(k)})\bar g(u^{j})]
\\
&=& \ds
\IE^{\rho^{J,l,\delta}}[\bar g(u^{(0)})^2]
+
2{1\over M}\sum_{k=0}^{M-1}\sum_{j=1}^{M-k}\cE^{\rho^{J,l,\delta},J,l}[\bar g(u^{(0)})\bar g(u^{(j)})]
\\
&=& \ds
\IE^{\rho^{J,l,\delta}}[\bar g(u^{(0)})^2]
+
2{1\over M}
\sum_{k=0}^{M-1}\sum_{j=1}^{M-k}\IE^{\rho^{J,l,\delta}}[\bar g(u^{(0)})\cE^{J,l}_{u^{(0)}}[\bar g(u^{(j)})]]
\\ 
&\le& \ds
\IE^{\rho^{J,l,\delta}}[\bar g(u^{(0)})^2]
\\ 
& & \ds
+2{1\over M}\sum_{k=0}^{M-1}\sup|\bar g|
\sum_{j=1}^{M-k}\IE^{\rho^{J,l,\delta}}[|\cE^{J,l}_{u^{(0)}}g(u^{(j)})-\IE^{\rho^{J,l,\delta}}[g]|]
\\
&\le& \ds
\IE^{\rho^{J,l,\delta}}[\bar g(u^{(0)})^2]
+
4{1\over M}\sum_{k=0}^{M-1}\sup|\bar g|^2 \sum_{j=1}^{M-k}(1-R)^j. 
\eeqas
In the last line we have used the estimate on the total variation
contraction from Proposition \ref{t:mc} noting, as in Remark \ref{r:uniform}, 
that because $\sup_{u\in U} \|P^{J,l}(u)\|_V$ is bounded uniformly 
with respect to the (discretization) parameters $J$ and $l$,
the constant $0 < R <1$ is independent of the parameters $J$ and $l$.
Since $\sup_{J,l} \IE^{\rho^{J,l,\delta}}[\bar g(u^{(0)})^2]$ 
is bounded independently of $J$ and of $l$,
we deduce that
\[
\sup_{J,l,M\in \IN} 
{M}\, \cE^{\rho^{J,l,\delta},J,l}
\Big[\Big|{\frac{1}{M}}\sum_{k=1}^M\bar g(u^{(k)})\Big|^2\Big]
<\infty 
\;.
\]
Next, we note that
\beqas
\cE^{\rho,J,l}\Big[\Big|\sum_{k=1}^M\bar g(u^{(k)})\Big|^2\Big]
&=&\ds
\int_U\cE^{J,l}_{u^{(0)}}\Big[\Big|\sum_{k=1}^M\bar g(u^{(k)})\Big|^2\Big]{d\rho(u^{(0)})}
\\
&=& \ds
\int_U\cE^{J,l}_{u^{(0)}}\Big[\Big|\sum_{k=1}^M\bar g(u^{(k)})\Big|^2\Big]\frac{d\rho}{d\rho^{J,l,\delta}}(u^{(0)})
d\rho^{J,l,\delta}(u^{(0)})
\\
&\le&\ds 
\cE^{\rho^{J,l,\delta},J,l}
\Big[\Big|\sum_{k=1}^M\bar g(u^{(k)})\Big|^2\Big]Z^{J,l}(\delta)\sup_{u\in U}
\exp(\Phi^{J,l}(u;\delta))
\;.
\eeqas
As $Z^{J,l}(\delta)\le 1$ and 
{$\sup_{u\in U}|\Phi^{J,l}(u;\delta)|$} 
are bounded uniformly with respect to $J$ and $l$, 
{we get the conclusion after using the bound from
Proposition \ref{prop:approxmu} on the
Hellinger distance between $\rho^{\delta}$ and 
$\rho^{J,l,\delta}.$}
\hfill$\Box$

We consider the case where 
$g(u)=\ell(P(\cdot,u))$ with $\ell\in V^*$ 
being a bounded linear functional on $V$. 
As there exists a constant $C>0$ such that for all
$J,l\in \mathbb{N}$, $u\in U$, there holds
\[
|\IE^{\rho^\delta}[\ell(P(\cdot, u))]-\IE^{\rho^\delta}[\ell(P^{J,l}(\cdot,u))]|
\le C(J^{-q}+2^{-l})
\] 
and 
\[
|\IE^{\rho^\delta}[\ell(P^{J,l}(\cdot,u))]-\IE^{\rho^{J,l,\delta}}[\ell(P^{J,l}(\cdot,u))]|
\le C(J^{-q}+2^{-l}),
\] 
we have 
\[
|\IE^{\rho^\delta}[\ell(P(\cdot,u))]-\IE^{\rho^{J,l,\delta}}[\ell(P^{J,l}(\cdot,u))]|
\le C(J^{-q}+2^{-l})\;.
\]
We therefore perform an MCMC algorithm to approximate 
$\IE^{\rho^{J,l,\delta}}[\ell(P^{J,l}(\cdot,u))]$. 
As $\ell(P^{J,l}(\cdot,u))$ and $\Phi^{J,l}(u;\delta)$ 
depend only on the finite set of coordinates $\{u_1,\ldots,u_J\}$ 
in expansion \eref{eq:a}, we perform the Metropolis-Hastings MCMC method on this set with 
proposals being drawn from the restriction of the prior measure $\rho$ to this finite set.   
\begin{proposition} \label{prop:fecost}
Let Assumptions \ref{assump:sumpsi} (i), (ii) and \ref{assump:FEcomplexity} hold,
and assume that $g(u)=\ell(P(\cdot, u))$ where $\ell$ is a bounded linear functional in $V^*$.
Then the approximate evaluation of the sample average 
${1\over M}\sum_{k=1}^M \ell(P^{J,l}(\cdot, u^{(k)}))$ 
by the Plain MCMC method with $M$ realizations of the chain, 
with Finite Element discretization in the domain $D$ 
at mesh level $l$ as described above, and with $J$-term truncated
coefficient representation \eref{eq:KtruncJ},
requires ${\mathcal O}(l^{d-1}2^{dl}MJ)$ floating point operations.
\end{proposition}
{\it Proof}\ \ 
The approximate computation of 
each of the $O(l^{d-1}2^{dl})$ non-zero entries 
of the stiffness matrix of \eref{eq:fe}, requires 
at most $O(J)$ operations to 
compute the coefficients $K^J$ at the quadrature points. 
\footnote{There is an implicit assumption here that the
basis functions can be evaluated at a point with ${\cal O}(1)$ cost.}
Therefore the cost of constructing the stiffness matrix is $O(l^{d-1}2^{dl}J)$. 
Since the condition number of this stiffness matrix
is assumed to be uniformly bounded,  
the conjugate gradient method for the approximate solution of
the linear system  
resulting from the Finite Element discretization
with an accuracy comparable to the order of the discretization
error requires at most $O(l^{d-1}2^{dl})$ floating point operations. 
The total cost for solving the approximated forward problem at 
each step of the Markov chain requires  at most
$O(l^{d-1}2^{dl}J)$ floating point operations. 
The numerical evaluation of 
$\ell(P^{J,l}(u^{(k)}))$ requires at most 
$O(2^{dl})$ floating point operations.
Since we generate $M$ draws of the chain, the assertion follows.
\hfill$\Box$
%
%
\begin{theorem} \label{t:fpe}
Let Assumptions \ref{assump:sumpsi} (i), (ii)
and \ref{assump:FEcomplexity} hold. 
For $g(u)=\ell(P(u))$ where $\ell$ is a bounded linear functional in $V^*$, 
with probability
$p_{\eN}(t)$ the conditional expectation $\IE^{\rho^\delta }g(u)$ 
can be approximated using $\eN$ degrees of freedom to 
approximate the forward PDE and $t^2 \eN^{2/d}$ MCMC steps 
(with a total of $t^2\eN^{1+2/d}$ degrees of freedom), 
incurring an error of $O(\eN^{-1/d})$, 
and  using at most 
$$c t^2\log(\eN)^{d-1}\eN^{1+(2+1/q)/d}$$ 
floating point operations for a positive constant $c$ where, for all $t$,
$$
\lim_{\eN \to \infty} p_{\eN}(t) \to \int_{-c't}^{c't}\frac{1}{\sqrt{2\pi}}\exp(-x^2/2)dx\;,
$$
for a positive constant $c'$ independent of $\eN$ and $t$. 

In mean square with respect to the measure ${\cal P}^{\rho,J,l}$, 
$\IE^{\rho^\delta}[g(u)]$ can be approximated with an error 
$O(N_{dof}^{-1/d})$, 
using at most $\eN^{1+2/d}$ number of degrees of freedom in total, 
and $O(\log(\eN)^{d-1}\eN^{1+(2+1/q)/d})$ floating point operations. 

\end{theorem}
{\it Proof}\ \ 
We first prove the probabilistic convergence result. 
We invoke the 
error estimate in Proposition \ref{prop:KJh}, 
and choose the parameters $M$, $J$ and $l$
so as to balance the bounds
$M^{-1/2}$, $J^{-q}$ and $2^{-l}$, 
taking into account
the fact that the coefficient of $M^{-\frac12}$
is only known through its asymptotic normality. 
We
select $J=2^{l/q}$ and $M=t^2N_{dof}^{2/d}$ where $t=c_3|\xi_M|$, with 
$N_{dof}$ denoting the number of degrees of freedom at each step  
being $N_{dof}=O(2^{dl})$; 
the constant $c_3$ and the random variable $\xi_M$ is as in Proposition \ref{prop:KJh}. 
Then the total number of floating point operations required 
as $l\rightarrow \infty$ is not larger than
$O(t^2 l^{d-1}2^{(d+2+1/q)l})$. 
We then arrive at the conclusion. 
The mean square convergence is proved in a similar manner.
\hfill$\Box$
%
\section{{Sparse gpc-}MCMC}
\label{sec:I2}
We again study computational complexity of the MCMC method 
defined by \eref{eq:MC} to sample the posterior probability 
distribution $\rho^\delta$. 
We adapt the approach in the previous section to use a computational method which effects
a reduction in computational cost by precomputing the
parametric dependence of the forward model, which enters
the likelihood. The main results are summarized in Theorem 
\ref{thm:improvedMCMC}. This method is introduced, and used in
practice, in the series of papers 
\cite{marzouk2007stochastic,marzouk2009stochastic,marzouk2009dimensionality}.
The major cost in MCMC methods 
is the repeated solution of the forward equations, 
with varying coefficients from the MCMC sampler
of $\rho^\delta$. 
The complexity of these repeated forward 
solves can be drastically reduced by 
precomputing an approximate, deterministic parametric representation
of the system's response which is valid for {\em all} 
possible realizations of $u\in U$. 
Specifically,
we precompute a sparse tensor finite element approximation
of the parametric, deterministic forward problem
by an approximate polynomial chaos representation of 
the solution's dependence on $u$ and by discretization 
of the forward solutions' spatial dependence from a 
multilevel hierarchy of Finite Element spaces in $D$.
This precomputed ``surrogate'' of the parametric response
map is then evaluated when running $M$ steps of the 
chain to estimate expectations with respect to the (approximate
due to the use of the surrogate) posterior measure.

As we shall show, this strategy is particularly effective 
if {the observations consist only of continuous}, 
linear functionals  of the system's response. 
In this case, only scalar coefficients
of the {forward map's} gpc expansion need 
to be stored and evaluated.
We use this to reduce the cost per step of the MCMC method. 
We continue to work under Assumption \ref{assump:sumpsi}(i),(ii)
and, furthermore, we add Assumption \ref{assump:sumpsi}(iii)
throughout the remainder of the paper.

\subsection{Approximation of the Posterior Measure}
For the solution $P_{N,{\mathcal L}}$ in Assumption \ref{ass:gpcplxity}, 
we define  
{the parametric, deterministic forward map based 
on the $N$-term truncated gpc expansion and Finite Element projected
surrogate forward map, 
as specified in Section \ref{sec:SPTApprox}
}
\be\label{eq:defGNcalL}
{\mathcal G}^{N,{\mathcal L}}(u)
=
(\cO_1(P_{N,{\mathcal L}}(\cdot,u)),\ldots,\cO_k(P_{N,{\mathcal L}}(\cdot,u))),
\ee
and the corresponding approximate Bayesian potential
\[
\Phi^{N,{\mathcal L}}(u;\delta)
=
{1\over 2}|\delta - {\mathcal G}^{N,{\mathcal L}}(u)|_{\Sigma}^2\;.
\]
The conditional measure $\rho^{N,{\mathcal L},\delta}$ 
on the measurable space $(U,\Theta)$ is defined as
\[
{d\rho^{N,{\mathcal L},\delta }\over d\rho}
\propto
\exp(-\Phi^{N,{\mathcal L}}(u;\delta))\;.
\]
We then have the following approximation result.
\begin{proposition}\label{prop:sparseHelldist}
Let Assumptions \ref{assump:sumpsi}(i) and \ref{ass:gpcplxity} hold.
Then there is a constant $C(\delta)$ which only depends on 
the data $\delta$ such that,
{
for every $N$ and $\mathcal L$ as in Assumption \ref{ass:gpcplxity},
}
\[
d_{\rm{Hell}}(\rho^\delta , \rho^{N,{\mathcal L},\delta})
\le 
C(\delta) 
\eN^{-\tau}
\;.
\]
\end{proposition}
{\it Proof}\ \ 
The proof for this proposition is similar to the
proof of Proposition \ref{prop:approxmu}, 
differing only in a few details; 
hence we highlight only the differences.
These are due to estimates on 
the forward error from Assumptions \ref{ass:gpcplxity} 
being valid only in the mean square sense 
whilst Proposition \ref{p:totapp} 
holds pointwise for $u \in U$. 
Nonetheless, at the point in the
estimation of $I_1$ and $I_2$ where the forward error estimate is
used, it is possible to use a mean square forward error estimate
instead of a pointwise forward error estimate. 
From Assumption \ref{ass:gpcplxity}, 
we deduce that there is a positive constant $c$ 
such that:
\[
\rho\{u:\ |\cG(u)-\cG^{N,\cL}(u)|>1\}\le cN_{dof}^{-2\tau}\;.
\]
As $\|P(u)\|_V$ is uniformly bounded for all $u$, there 
is a constant $c_1(\delta)>0$ such that 
$|\delta-\cG(u)|_\Sigma < c_1(\delta)$. 
Choose a constant  $c_2(\delta)>0$  sufficiently large. 
If $|\delta-\cG^{N,\cL}(u)|_\Sigma>c_2(\delta)$, 
then
\[
|\cG^{N,\cL}(u)-\cG(u)|_\Sigma
\ge 
|\delta-\cG^{N,\cL}(u)|_\Sigma-|\delta-\cG(u)|_\Sigma
> c_2(\delta)-c_1(\delta)>1\;.
\]
Let $U_1\subset U$ be the set of $u\in U$ 
such that $|\delta-\cG^{N,\cL}(u)|_\Sigma>c_2(\delta)$. 
We have that
$ \rho(U_1)\le cN_{dof}^{-2\tau}$.
Thus,
\[
{1\over Z(\delta)}
\int_{U_1}
\Big|\exp\Big(-{1\over 2}\Phi(u;\delta)\Big)-\exp\Big(-{1\over 2}\Phi^{N,\cL}(u;\delta)\Big)\Big|d\rho(u)
\le C(\delta)N_{dof}^{-2\tau}
\;.
\]
When $u\notin U_1$, $|\delta-\cG^{N,\cL}(u)|_\Sigma\le c_2(\delta)$ 
so there is a constant 
$c_3(\delta)$ so that $|\cG^{N,\cL}(u)|\le c_3(\delta)$.
An argument similar to that for \eref{eq:LipschitzPhi} shows that 
\beqas
\Big|\exp\left(-{1\over2}\Phi(u;\delta)\right)-\exp\left(-{1\over2}\Phi^{N,{\mathcal L}}(u;\delta)\right)\Big|\le\qquad\qquad \\
C(2|\delta|+|\cG(u)|+|\cG^{N,{\mathcal L}}(u)|)|\cG(u)-\cG^{N,{\mathcal L}}(u)|.
\eeqas
Therefore
\beqas
I_1&=&{1\over Z(\delta)}\int_U\Big|\exp\Big(-{1\over 2}\Phi(u;\delta)\Big)-\exp\Big(-{1\over 2}\Phi^{N,\cL}(u;\delta)\Big)
\Big|^2d\rho(u)
\\
&\le& C(\delta)N_{dof}^{-2\tau}
+
\\
&& c\int_U(2|\delta|+|\cG(u)|+c_3(\delta))^2|\cG(u)-\cG^{N,{\mathcal L}}(u)|^2d\rho(u)
\\
&\le& C(\delta)N_{dof}^{-2\tau}+C(\delta)\int_U\|P(\cdot,u)-P_{N,{\mathcal L}}(\cdot,u)\|_V^2d\rho(u)
\\
&\le&C(\delta)N_{dof}^{-2\tau}
\;.
\eeqas
To show that $I_2 < C(\delta)\eN^{-2\tau}$ we still need to verify that 
\[
Z^{N,{\mathcal L}}(\delta)=\int_U\exp(-\Phi^{N,{\mathcal L}}(u;\delta))d\rho(u)
\]
is uniformly bounded from below by a positive bound for all $N$ and $\mathcal L$. 
As $P_{N,{\mathcal L}}$ is uniformly bounded in $L^2(U,\rho;V)$, 
\[
\int_U|\cG^{N,{\mathcal L}}(u)|d\rho(u)\le c\int_U\|P_{N,{\mathcal L}}(u)\|_Vd\rho(u)\le c.
\]
Fixing $r>0$ sufficiently large, 
the $\rho$ measure of the set $u\in U$ such that $|\cG^{N,{\mathcal L}}(u)|>r$ 
is bounded by $c/r$. 
Therefore the measure of the set of $u\in U$ such that $|\cG^{N,{\mathcal L}}(u)|\le r$ is 
bounded from below by $1-c/r$. 
Thus we have proved that 
for every realization $\delta$ of the data,
there exists $C(\delta)>0$ such that
\[
Z^{N,{\mathcal L}}(\delta)
\ge\int_U\exp(-{1\over2}(|\delta|_\Sigma
+
|\cG^{N,{\mathcal L}}(u)|_\Sigma)^2)d\rho(u)
>
C(\delta)
>0
\;.
\]   
\hfill$\Box$

%
Let $(u^{(k)})_k$ be the Markov chain generated by the 
sampling process \eref{eq:MC} with the 
acceptance probability being replaced by
\be
\alpha^{N,\cL}(u,v) = 1\wedge\exp(\Phi^{N,\cL}(u,\delta)-\Phi^{N,\cL}(v,\delta)) \;.
\label{eq:alphaNcL}
\ee
We denote by
\[
E_M^{\rho^{N,\cL,\delta}}[g]={1\over M}\sum_{k=1}^Mg(u^{(k)})\;.
\]
{
We then have the following error for the gpc-FE surrogate based 
MCMC method (for notation we refer to Section \ref{ssec:OvNotat})
}
\begin{proposition}\label{prop:Ephiu} 
Let Assumptions \ref{assump:sumpsi}(i) and \ref{ass:gpcplxity} hold and
let $g$ be a bounded continuous function from $U$ to $\mathbb{R}$. 
Then
\be
|\IE^{\rho^\delta}[g]-E_M^{\rho^{N,\cL,\delta}}[g]|\le c_6M^{-1/2}+c_7\eN^{-\tau},
\label{eq:Ephiu}
\ee
$\cP^{\rho^{N,\cL,\delta},N,\cL}$ almost surely,
where $c_6\le c_8|\xi_M|$ where $\xi_M$ is a random variable which 
converges weakly as $M\to\infty$ to $\xi\sim N(0,1)$; 
the constants $c_7$ and $c_8$ are deterministic and do not depend on $M$, $N$ and $\eN$.

There exists a deterministic positive constant $c_9$  
such that the gpc-MCMC converges in the mean square with the same rate of convergence
\[
\Big(\cE^{\rho,N,\cL}\Big[\left|\IE^{\rho^\delta}[g]-E_M^{\rho^{N,\cL,\delta}}[g]\right|^2\Big]\Big)^{1/2}
\le 
c_9(M^{-1/2}+N_{dof}^{-\tau})
\;.
\]
\end{proposition}
{\it Proof}\ \ 
Using \eref{eq:alphaNcL}, a random draw from $\rho$ 
has probability larger than $\exp(-\Phi^{N,\cL}(v;\delta))$ of being accepted. 
Therefore the transition kernel of the Markov chain generated by \eref{eq:MC} 
with the acceptance probability \eref{eq:alphaNcL} satisfies
\[
p(u,A)\ge \int_A\exp(-\Phi^{N,\cL}(v;\delta))d\rho(v).
\]
Using Theorem 16.2.4 of \cite{MT93}, 
we deduce that the $n$th iteration of the transition kernel satisfies
\[
\|p^n(u,\cdot)-\rho^{N,\cL,\delta}\|_{\rm{TV}}
\le 2
\left(1-\int_U\exp(-\Phi^{N,\cL}(v;\delta))d\rho(v)\right)^n
\;.
\]
From the proof of Proposition \ref{prop:sparseHelldist}, we have
\[
\int_U\exp(-\Phi^{N,\cL}(v;\delta))d\rho(v)
\ge 
\exp(-c_2(\delta)^2/2)+cN_{dof}^{-2\tau}\;.
\]
Thus, we can choose a constant $R<1$ independent of the approximating parameters $N$ and $\mathcal L$
so that for all $n\in \IN$ holds
\[
\|p^n(u,\cdot)-\rho^{N,\cL,\delta}\|_{\rm{TV}}\le 2(1-R)^n
\;.
\]
In a similar manner as for Proposition \ref{prop:KJh}, we deduce the probabilistic bound. 
For the mean square bound, similar to the proof of Proposition \ref{prop:KJh}, we have
\[
\cE^{\rho^{N,\cL},N,\cL}\Big[\left|\IE^{\rho^\delta}[g]-E_M^{\rho^{N,\cL,\delta}}[g]\right|^2\Big]
\le 
C(M^{-1/2}+N_{dof}^{-\tau})^2\;.
\]
Let $U':=\{u\in U :\ |\cG^{N,\cL}(u)-\cG(u)|>1\}$. 
We deduce that 
there exists a constant $c>0$ independent of $\cL$, $N_{dof}$, $N$ such that
$\rho(U')\le c N_{dof}^{-2\tau}$ and such that we may estimate
\[
\begin{array}{rl}
& \ds
\cE^{\rho,N,\cL}\Big[\left|\IE^{\rho^\delta}[g]-E_M^{\rho^{N,\cL,\delta}}[g]\right|^2\Big]
\\
=&\ds
\int_{U'}\cE_{u^{(0)}}^{N,\cL}\Big[\left|\IE^{\rho^\delta}[g]-E_M^{\rho^{N,\cL,\delta}}[g]\right|^2\Big]{d\rho(u^{(0)})}
\\
& \ds \qquad 
+ 
\int_{U\setminus U'}
\cE_{u^{(0)}}^{N,\cL}\Big[\left|\IE^{\rho^\delta}[g]-E_M^{\rho^{N,\cL,\delta}}[g]\right|^2\Big] d\rho(u^{(0)})
\\
\le& \ds
C N_{dof}^{-2\tau}+\int_{U\setminus U'}\cE_{u^{(0)}}^{N,\cL}\Big[\left|\IE^{\rho^\delta}[g]-E_M^{\rho^{N,\cL,\delta}}[g]\right|^2\Big]d\rho(u^{(0)})
\\
\le& \ds 
C N_{dof}^{-2\tau}+
\\
& \ds
\int_{U\setminus U'}
\cE_{u^{(0)}}^{N,\cL}
\Big[\left|\IE^{\rho^\delta}[g]-E_M^{\rho^{N,\cL,\delta}}[g]\right|^2\Big]
Z^{N,\cL}(\delta)\exp(\Phi^{N,\cL}(u;\delta))
d\rho^{N,\cL,\delta}(u^{(0)})
\;.
\end{array}
\]
On $U\setminus U'$, $\sup_{u\in U}|\cG^{N,\cL}(u)|$ is uniformly bounded 
with respect to all $N$ and $\cL$. 
From this, we arrive at the conclusion.
\hfill$\Box$
%
%
{
\begin{remark}\label{rem:5.9}
In Proposition \ref{prop:Ephiu}, $g$ is assumed to 
be a bounded continuous function from $U$ to $\IR$. 
The {gpc-FE surrogate accelereated} 
MCMC is of particular interest in the case
where g is given by $\ell \circ P$ where 
$\ell$ is a bounded, linear functional on $V$, i.e. $\ell \in V^*$. 
From Assumption \ref{ass:gpcplxity} and the fact that 
\[
{d\rho^\delta\over d\rho}{(u)}
=
{1\over Z(\delta)}\exp(-\Phi(u;\delta))
\;,
\]
we deduce that
\[
|\IE^{\rho^\delta}[\ell(P(\cdot,u))]-\IE^{\rho^\delta}[\ell(P_{N,\cL}(\cdot,u))]|
\le 
C(\delta,\ell) N_{dof}^{-\tau}
\;.
\]
On the other hand, from Proposition \ref{prop:sparseHelldist}, 
we have ({cf. \cite[Eq. (2.7)]{Cotteretal}})
\[
|\IE^{\rho^\delta}[\ell(P_{N,\cL}(\cdot, u))]-\IE^{\rho^{N,\cL,\delta}}[\ell(P_{N,\cL}(\cdot,u))]|
\le C(\delta,\ell) N_{dof}^{-\tau}\;.
\]
Therefore,  by the triangle inequality,
\[
|\IE^{\rho^\delta}[\ell(P(\cdot,u))]-\IE^{\rho^{N,\cL,\delta}}[\ell(P_{N,\cL}(\cdot,u))]|
\le C(\delta,\ell)N_{dof}^{-\tau}
\;.
\]
{
We wish to
approximate $\IE^{\rho^{N,\cL,\delta}}[\ell(P_{N,\cL}(\cdot,u))]$
with a Markov Chain-Monte Carlo algorithm.
}
In doing so, the following difficulty may arise:
although $\ell(P(\cdot,u))$ is uniformly bounded with respect to $u\in U$, 
$\sup_{u\in U} \ell(P_{N,\cL}(\cdot,u))$ may not be uniformly bounded 
with respect to $N$ and $\cL$. 
However, we can still apply Proposition \ref{prop:Ephiu}  
by using a cut-off argument:
to this end, we define the continuous bounded function 
$\tilde g(u):U\to\IR$ by truncation, 
i.e. 
\[
\fl\tilde g(u) 
:= 
\left\{ 
\begin{array}{ccr}
\ds
\ell(P_{N,\cL}(\cdot,u)) &\mbox{ if }& |\ell(P_{N,\cL}(\cdot,u))|
\le \sup_{u\in U}{|\ell(P(\cdot,u))|}+1 \;,
\\
\ds
\sup_{u\in U} |\ell(P(\cdot,u))|+1 &\mbox{ if }& \ell(P_{N,\cL}(\cdot,u))
> 
\sup_{u\in U}{|\ell(P(\cdot,u))|}+1 \;,
\\
\ds
-\sup_{u\in U}|\ell(P(\cdot,u))|-1 &\mbox{ if }& \ell(P_{N,\cL}(\cdot,u)) 
< 
-\sup_{u\in U}{|\ell(P(\cdot,u))|}-1
\;.
\end{array}
\right.
\]
Define $U':=\{ u\in U: |\ell(P(\cdot,u))-\ell(P_{N,\cL}(\cdot,u))|>1\}$. 
From Assumption \ref{ass:gpcplxity}, 
we find that $\rho(U')<C(\ell,\delta)N_{dof}^{-2\tau}$.
It follows then that there exists a constant $c>0$ depending
on the data $\delta$, but independent of $N$ and of $\cL$
such that
\beqas
|\IE^{\rho^{N,\cL,\delta}}[\ell(P_{N,\cL}(\cdot,u))-\tilde g(u)]|
\le 
\int_{U'}|\ell(P_{N,\cL}(\cdot,u))-\tilde g(u)|d\rho^{N,\cL,\delta}(u)
\\
\le C(\delta)\int_UI_{U'}(u)(|\ell(P_{N,\cL}(\cdot,u))|+c)d\rho(u)
\\
\le C(\delta)\rho(U')^{1/2}(\|\ell(P_{N,\cL}(\cdot,u))\|_{L^2(U,\rho;\IR)}+c)
\le C(\delta)N_{dof}^{-\tau}
\;.
\eeqas
Therefore, we may run the MCMC algorithm on $\IE^{\rho^{N,\cL,\delta}}[\tilde g(u)]$. 
\end{remark} 
}
At each step of the MCMC algorithm, we need to compute 
$\ell(P_{N,\cL}(\cdot,u^{(k)}))$ 
which, 
{\em for linear functionals} $\ell(\cdot)$, is equal to
$\sum_{\nu\in \Lambda}\ell(P_{\nu,\cL})L_\nu(\cdot,u^{(k)})$.
Because the parametric solution of the elliptic problem can
be precomputed before the MCMC is run, and then needs
only to be {\em evaluated} at each state of the MCMC method,
significant savings can be obtained. We illustrate this,
using the ideas of the previous Remark \ref{rem:5.9}, 
to guide the choice of test functions.
\begin{proposition}\label{prop:sparsecomplexity}
Let Assumptions \ref{assump:sumpsi}(i) and \ref{ass:gpcplxity} hold
and let $g(u)=\ell(P(\cdot,u))$ where $\ell\in V^*$. 
Then the total number of floating point 
operations required for performing $M$ 
steps in the Metropolis-Hastings method as $N,M \rightarrow \infty$
is bounded by $O(\eN^\alpha (\log \eN)^\beta+MN\log N)$.
\end{proposition}
{\it Proof}\ \ 
By Assumption \ref{ass:gpcplxity} and with the notation as in that assumption,
the cost of solving one instance of problem \eref{eq:parameterfe} 
is bounded by 
$O(\eN^\alpha (\log \eN)^\beta)$.
At each MCMC step, we need to evaluate
the observation functionals
\beq \label{eq:gpcOutput}
\cO_i(P_{N,{\mathcal L}}(u^{(k)}))
=
\sum_{\nu\in\Lambda_N}\cO_i(P_{\nu,\mathcal{L}})L_\nu(u^{(k)})
\;.
\eeq
{
We note in passing that the storage of the parametric 
gpc-type representation of the forward map \eref{eq:gpcOutput}
requires only one real per gpc mode, {\em provided that
only functionals of the forward solution are of interest}.
We now estimate the complexity of computing one draw 
of the forward map \eref{eq:gpcOutput}.
}
For $\nu\in \cF$, 
each multivariate Legendre polynomial 
$L_\nu(u^{k})$ can be evaluated with $O(|\nu|)$ float point
operations.
As $|\nu|=O(\log N)$, computing the observation
functionals $\cO_i(P_{N,{\mathcal L}})$ requires
$O(N\log N)$ floating point operations. 
Thus we need $O(\eN^\alpha (\log \eN)^\beta + MN\log N)$ floating point operations to perform $M$ steps of the Metropolis-Hastings method
with sampling of the surrogate, sparse gpc-Finite Element
representation of the forward map.

\hfill$\Box$
\begin{theorem}\label{thm:improvedMCMC}
Let Assumptions \ref{assump:sumpsi}(i) and \ref{ass:gpcplxity} hold.
For $g(u)=\ell(P(\cdot,u))$ with given $\ell\in V^*$, 
with probability $p_{\eN}(t)$ the conditional expectation 
$\IE^{\rho^\delta}[g(u)]$ can be approximated with $\eN$ degrees
of freedom, incurring an error of $O(\eN^{-\tau})$ using
at most 
\[
c\eN^\alpha(\log \eN)^\beta+ct^2\eN^{2\tau+\tau/\sigma}\log(\eN)
\]
many floating point operations,  where
\[
\lim_{\eN\to\infty}p_{\eN}(t)\to\int_{-c't}^{c't}{1\over\sqrt{2\pi}}\exp(-x^2/2)dx,
\]
for some  constants $c,c'$ independent of $\eN$.

In the mean square with respect to the measure ${\cal P}^{\rho,N,{\cal L}}$, $\IE^{\rho^\delta}[g(u)]$ 
can be approximated with $N_{dof}$ degrees of freedom, with an error $N_{dof}^{-\tau}$ 
using at most 
\[
O(\eN^\alpha(\log\eN)^\beta + \eN^{2\tau+\tau/\sigma}\log(\eN))
\]
floating point operations. 
\end{theorem}
{\it Proof}\ \ 
{
We relate the number of MCMC realizations $M$ with the total
number of degrees of freedom $N_{dof}$ by equating the terms in the 
error bound \eref{eq:Ephiu}.
To this end,
}
we choose 
{
$M=t^2\eN^{2\tau}$ where $t=c_8|\xi_M|$; 
the constant $c_8$ and the random variable $\xi_M$ are as in Proposition \ref{prop:Ephiu} . 
With $N =  O(\eN^{\tau/\sigma}) $, 
the number of floating point operations required in Proposition \ref{prop:sparsecomplexity} 
is bounded by
$$
O\left(\eN^\alpha(\log \eN)^\beta + t^2\eN^{2\tau+\tau/\sigma}\log \eN\right)
\;.
$$
As $\xi_M$ converges weakly to the normal Gaussian variable, 
we deduce the limit for the probability density $p_{\eN}(t)$ of the random variable $t$. 
}
The proof for the mean square approximation is similar. 
\hfill$\Box$

\begin{remark} 
In studying complexity of the Plain MCMC method in Theorem \ref{t:fpe}, 
the discretized parametric PDE \eref{eq:prob} 
is to be solved once at every step of the MCMC process, 
using $\eN$ degrees of freedom, with $O(\eN^{2/d})$ 
steps required (the multiplying constant depends 
on a random variable when we consider the realization-wise error). 
Ignoring $\log$ factors, 
the error resulting from discretization and
running the MCMC on the discretized PDE
can be bounded in terms of the total 
number of floating point operations $N_{fp}$ 
by $O(N_{fp}^{-1/(d+2+1/q)})$. 
In Theorem \ref{thm:improvedMCMC} the discretized forward PDE 
is solved for every realization before running the MCMC process. 
The rate of convergence of the MCMC process in terms of the 
total number of floating point operations used 
is $O(N_{fp}^{-{\min(\tau/\alpha,1/(2+1/\sigma))}})$. 
This can be significantly smaller than the rate of 
convergence in Theorem \ref{t:fpe} when $\alpha$ is 
close to 1.  To see this we need to dig into the manner 
in which Assumptions \ref{ass:gpcplxity} are typically verified
which, in turn, requires Assumptions \ref{assump:sumpsi}(iii),(iv)
(which together imply Assumptions \ref{assump:sumpsi}(ii), see
discussion in Appendix E).
For example, with the decay rate of $\|\psi_j\|_\infty$ in Assumption \ref{assump:sumpsi}(iii), 
the summability constant $p$ in Assumption \ref{ass:gpclp} can be any constant that is greater than $1/s$. 
Therefore, the gpc approximation rate 
$\sigma$ in Proposition \ref{prop:Stechkin} 
can be chosen as any positive constant smaller than $s-1/2$. 
On the other hand, 
the $J$-term approximation rate $q$ 
of the unknown input $K$ in Assumption \ref{assump:sumpsi}(ii) is bounded by $s-1$.  
As 
\[
2+{1\over s-1/2}<d+2+{1\over s-1},
\]
we therefore can choose $\sigma$ so that
\[
{1\over 2+1/\sigma}>{1\over {d+2+1/q}}\;.
\]
As shown in \cite{CDS1}, when $(\|P_\nu\|_{H^2(D)})_\nu\in\ell^p(\cF)$, 
$\tau$ can be chosen as $1/d$. 
Thus, when $\alpha$ is sufficiently close to 1, the complexity of the sparse 
gpc-MCMC approach is superior 
to that of the Plain MCMC approach in the previous section.
\end{remark}
\section{Multilevel MCMC}
\label{sec:M}
In the preceding section
we showed that complexity reduction 
is possible in the Plain MCMC sampling of 
the posterior measure $\rho^\delta$
provided that all samples are determined 
from {\em one} precomputed, ``surrogate'' 
sparse tensor gpc-representation of the 
forward map of the discretization of the 
parametric, deterministic problem \eref{eq:prob}.
Specifically, 
we proved that when the forward map $\cG(u)$
consists of {\em continuous, linear functionals $\cO_i(\cdot)$ 
on the forward solution $U\ni u\mapsto P(\cdot,u) \in V$}, 
 {and when the functional whose posterior expectation we seek
is also linear on this space},
this gpc-type approximation of the solution can reduce 
the complexity required.
Lower efficiency results if, 
for example, the rate of convergence 
of the sparse tensor finite element
solution in \eref{eq:parameterfe} is 
moderate in terms of
the total number of degrees of freedom, 
and/or if the complexity grows superlinearly with respect 
to the number of degrees of freedom.
Furthermore, although an increasing number of efficient 
algorithms for the computation 
of approximate responses of the forward problem on the 
entire parameter space $U$ are available (e.g. \cite{BAS,ScCJG,CJGdiss,b:09,CCDS11})
and therefore gpc-surrogates for the MCMC are available, 
many systems of engineering interest may not admit 
gpc-based representations of the parametric forward maps.
Finding other, non-gpc based, methods for reducing 
the complexity of Plain MCMC sampling under 
$\rho^\delta$ is therefore of interest. 
We do this by using ideas from Multilevel Monte Carlo.
The resulting complexity is summarized in Theorem \ref{thm:MLMCMCerr}.
We give sufficient conditions 
on the approximation methods and on the {basis functions $\psi_j$ 
appearing in \eref{eq:a}} such that 
complexity reduction is possible by performing a 
multilevel sampling procedure where 
a number of samples depending on the discretization 
parameters is used for problem \eref{eq:prob}.

\subsection{Derivation of the MLMCMC}
\label{ssec:DerivMLMCMC}
For given, fixed $\ell\in V^*$, a bounded linear functional on $V$,
we aim at estimating $\IE^{\rho^\delta}[\ell(P(\cdot,u))]$ 
where $P$ is the solution of problem \eref{eq:prob}. 
For each level $l = 1,2 ,\ldots,L$, 
we assume that problem \eref{eq:prob} is discretized with  
the truncation of the \KL expansion after $J$ terms
with $J=J_l$ as defined in \eref{eq:KtruncJ}
and with a finite element discretization meshwidth $h_l$.
The multilevel FE-discretization of the forward problem
\eref{eq:prob} and the truncation \eref{eq:KtruncJ} 
induces a corresponding hierarchy of approximations 
$\rho^{J_{l},l,\delta}$ of the posterior measure $\rho^\delta$.

Following \cite{G08,BSZ10,MS10,Brinkman} the MLMCMC will be based
on sampling a telescopic expansion of the discretization
error with a level-dependent sample size. 
Contrary to
\cite{G08,BSZ10,MS10}, however, we introduce 
now two multilevel discretization hierarchies, 
one for 
the parametric forward solutions $\{P^{J_l,l}\}_{l=0}^L$ 
(where the level corresponds to mode truncation of the coefficient
 and to discretization for the approximate solution
 of the parametric forward problem) and 
a second one $\{\rho(J_{l'},l',\delta)\}_{l'=0}^L$
for the posterior measure.
The `usual'  
telescoping argument as in \cite{G08} or in 
\cite{BSZ10,MS10} 
together with the single-level error bound in 
Proposition \ref{prop:KJh} alone 
does not allow for obtaining a convergence for the present problem.

As in Section \ref{sec:S}, we work under Assumptions 
\ref{assump:sumpsi} (i), (ii) and \ref{assump:FEcomplexity}.
We recall the sequence of discretization levels 
in the FE discretizations in the physical domain $D$,
as in Assumption \ref{assump:FEcomplexity},
and the input truncation dimension $J$ in 
Assumption \ref{assump:sumpsi}(ii).
We then derive the Multilevel MCMC-FEM as follows. 
First, we note that, by Assumption \ref{assump:FEcomplexity}, using
the uniform boundedness of the $P^J$ in $W$ as in the preceding
section, there exists $C>0$ independent of $L$ such that
\be
\fl|\IE^{\rho^\delta}[\ell(P(\cdot,u))]-\IE^{\rho^\delta}[\ell(P^{J_L,L}(\cdot,u))]|
\le 
C\sup_{u\in U}\|P(\cdot,u)-P^{J_L,L}(\cdot,u)\|_V
\le 
C2^{-L} 
\;.
\label{eq:1}
\ee
With the convention that {expectation with respect to 
$\rho^{J_{-1},-1,\delta}$ denotes integration with respect to the
measure which assigns zero mass to all subsets of the probability
space}, 
we write
\begin{eqnarray}
\IE^{\rho^{J_L,L,\delta}}[\ell(P(\cdot,u))]
 &= \ds
\sum_{l = 0}^L  
\left( \IE^{\rho^{J_l ,l,\delta}}[\ell(P(\cdot,u))] - \IE^{\rho^{J_{l-1},l-1,\delta}}[\ell(P(\cdot,u))] 
\right)
\\
&= \ds \sum_{l = 0}^L
\left(\IE^{\rho^{J_l ,l,\delta}}-\IE^{\rho^{J_{l-1},l-1,\delta}}\right)[\ell(P(\cdot,u))]
\;.
\label{eq:1z}
\end{eqnarray}
Analogously we may write, {for any $L' \le L$} 
(omitting the arguments of $P$ and its approximations for 
brevity of notation)
\begin{equation}
\label{eq:2z}
\IE^{\rho^{J_L,L,\delta}}[\ell(P^{J_{L'},L'})]
=
\ds \sum_{l = 0}^L
\left(\IE^{\rho^{J_l ,l,\delta}}-\IE^{\rho^{J_{l-1},l-1,\delta}}\right)[\ell(P^{J_{L'},L'})].
\end{equation}
With the convention that $P^{J_{-1},-1} := 0$
we have for each $l$ and $L'$
\begin{equation} \label{eq:PLtelescop}
\fl\left(\IE^{\rho^{J_l,l,\delta}}-\IE^{\rho^{J_{l-1},l-1,\delta}}\right) [\ell(P^{J_{L'},L'})] 
=
\sum_{l' = 0}^{L'} 
\left(\IE^{\rho^{J_l,l,\delta}}-\IE^{\rho^{J_{l-1},l-1,\delta}}\right)[\ell(P^{J_{l'},l'}) - \ell(P^{J_{l'-1},l'-1})]
\;.
\end{equation}
Finally, a computable
multilevel approximation will be obtained
on running, for each level $l=0,1,...,L$
of truncation resp. Galerkin projection,
simultaneously a suitable number 
of realizations of a Markov chain ${\mathcal C}_l$
to approximate the expecations
\be
\label{eq:1zz}
\begin{array}{rl}
& \displaystyle
\sum_{l = 0}^L
\left(\IE^{\rho^{J_l ,l,\delta}}-\IE^{\rho^{J_{l-1},l-1,\delta}}\right)[\ell(P^{J_{L'},L'})]
\\
=& \displaystyle
\sum_{l=0}^L\sum_{l' = 0}^{L'}
\left(\IE^{\rho^{J_l,l,\delta}}-\IE^{\rho^{J_{l-1},l-1,\delta}}\right)
[\ell(P^{J_{l'},l'}) - \ell(P^{J_{l'-1},l'-1})]
\end{array}
\ee
by sample averages of $M_{ll'}$ many realizations, upon
choosing $L'(l)$ judiciously.
To derive a computable MLMCMC estimator 
we observe that, for any {measurable} functions
$Q:U\to \IR$ which is integrable with respect to 
the approximate posterior measures $\rho^{J_l,l,\delta}$,
\beqas
&&\fl\left(\IE^{\rho^{J_l,l,\delta}}-\IE^{\rho^{J_{l-1},l-1,\delta}}\right)[Q]
\\
&&\fl={1\over Z^{J_l,l}}\int_U\exp(-\Phi^{J_l,l}(u;\delta))Q(u)d\rho(u)-{1\over Z^{J_{l-1},l-1}}\int_U\exp(-\Phi^{J_{l-1},l-1}(u;\delta))Q(u)d\rho(u)
\\
&&\fl={1\over Z^{J_l,l}}\int_U\exp(-\Phi^{J_l,l}(u;\delta))\left(1-\exp(\Phi^{J_l,l}(u;\delta)-\Phi^{J_{l-1},l-1}(u;\delta))\right)Q(u)d\rho(u)
\\
&&\fl\qquad\qquad+\left({Z^{J_{l-1},l-1}\over Z^{J_l,l}}-1\right){1\over Z^{J_{l-1},l-1}}\int\exp(-\Phi^{J_{l-1},l-1}(u;\delta))Q(u)d\rho(u)
\;.
\eeqas
We note further that
$$
{Z^{J_{l-1},l-1}\over Z^{J_l,l}}-1
=
{1\over Z^{J_l,l}}
\int_U
\left(\exp(\Phi^{J_l,l}(u;\delta)-\Phi^{J_{l-1},l-1}(u;\delta))-1\right)\exp(-\Phi^{J_l,l}(u;\delta))
d\rho(u)
\;.
$$
Thus an approximation for 
$Z^{J_{l-1},l-1}/ Z^{J_l,l}-1$ can be found by 
running the MCMC with respect to 
the approximate posterior
$\rho^{J_l,l,\delta}$ to sample 
the potential difference
$\exp(\Phi^{J_l,l}(u;\delta)-\Phi^{J_{l-1},l-1}(u;\delta))-1$.
We define the Multilevel Markov Chain Monte Carlo estimator 
$E^{MLMCMC}_L[\ell(P)]$ of $\IE^{\rho^{{\delta}}}[\ell(P)]$
as
$$
\begin{array}{l}
\displaystyle
E_L^{MLMCMC}[\ell(P)]
=
\\ \\
\displaystyle
\sum_{l=0}^L\sum_{l'=0}^{L'(l)}E_{M_{ll'}}^{\rho^{J_l,l,\delta}}
\left[\Big(1-\exp(\Phi^{J_l,l}(u;\delta)-\Phi^{J_{l-1},l-1}(u;\delta))
      \Big)(\ell(P^{J_{l'},l'})-\ell(P^{J_{l'-1},{l'-1}}))\right] 
\\ \\
\displaystyle
+ 
E_{M_{ll'}}^{\rho^{J_l,l,\delta}}
\left[\exp(\Phi^{J_l,l}(u;\delta)-\Phi^{J_{l-1},l-1}(u;\delta))-1\right]\cdot
E_{M_{ll'}}^{\rho^{J_{l-1},l-1,\delta}}\left[\ell(P^{J_{l'},l'})-\ell(P^{J_{l'-1},l'-1})
\right]
\;.
\end{array}
$$
Up to this point, 
the choice of the index $L'(l)$ 
and the sample sizes $M_{ll'}$ 
are still at our disposal.
{Choices for them will be made} based on an error-versus-work 
analysis of this estimator {that we now present}.
\subsection{Error Analysis}
\label{ssec:MLMCMCError}
{To perform} the error analysis of the MLMCMC approximation
we decompose the error into three terms as follows.
\begin{proposition}
\label{prop:MLMCMCerr}
There holds
\beq\label{eq:MLMCMCerr}
\IE^{\rho^{{\delta}}}[\ell(P)] - E^{MLMCMC}_L[\ell(P)]
= I_L + II_L + III_L
\eeq
where
$$
I_L := \IE^{\rho^{{\delta}}}[\ell(P)] - \IE^{\rho^{J_L,L,\delta}}[\ell(P)],
\quad
II_L=\sum_{l=0}^L(\IE^{\rho^{J_l,l,\delta}}-\IE^{\rho^{J_{l-1},l-1,\delta}})[\ell(P)-\ell(P^{J_{L'(l)},L'(l)})] 
$$
and
\beqas
\fl III_L 
:= 
\sum_{l=0}^L \sum_{l'=0}^{L'(l)}\bigg\{\IE^{\rho^{J_l,l,\delta}}\left[(1-\exp(\Phi^{J_l,l}(u;\delta)-\Phi^{J_{l-1},l-1}(u;\delta)))](\ell(P^{J_{l'},l'})-\ell(P^{J_{l'-1},l'-1}))\right]\\
-E_{M_{ll'}}^{\rho^{J_l,l,\delta}}\left[(1-\exp(\Phi^{J_l,l}(u;\delta)-\Phi^{J_{l-1},l-1}(u;\delta)))](\ell(P^{J_{l'},l'})-\ell(P^{J_{l'-1},l'-1}))\right]\bigg\}\\\
+\bigg\{\IE^{\rho^{J_l,l,\delta}}[\exp(\Phi^{J_l,l}(u;\delta)-\Phi^{J_{l-1},l-1}(u;\delta))-1]\cdot\IE^{\rho^{J_{l-1},l-1,\delta}}[\ell(P^{J_{l'},l'})-\ell(P^{J_{l'-1},l'-1})]\\
-E_{M_{ll'}}^{\rho^{J_l,l,\delta}}[\exp(\Phi^{J_l,l}(u;\delta)-\Phi^{J_{l-1},l-1}(u;\delta))-1]\cdot E_{M_{ll'}}^{\rho^{J_{l-1},l-1,\delta}}[\ell(P^{J_{l'},l'})-\ell(P^{J_{l'-1},l'-1})]\bigg\}
\eeqas
respectively.

%
\end{proposition}
{\it Proof}
From equation \eref{eq:1z} we have
\begin{equation}
\fl\IE^{\rho^\delta}[\ell(P)] - \IE^{\rho^{J_L,L,\delta}}[\ell(P)]
 = \ds
\IE^{\rho^\delta}[\ell(P)] 
-
\sum_{l = 0}^L  
\left( \IE^{\rho^{J_l ,l,\delta}}[\ell(P)] - \IE^{\rho^{J_{l-1},l-1,\delta}}[\ell(P)] \right)
\;.
\label{eq:11z}
\end{equation}
from which it follows that
\begin{eqnarray*}
\fl\IE^{\rho^{\delta}}[\ell(P)] - \IE^{\rho^{J_L,L,\delta}}[\ell(P)]
 &= \ds
\IE^{\rho^{\delta}}[\ell(P)]\\ 
&-
\sum_{l = 0}^L  
\left( \IE^{\rho^{J_l ,l,\delta}}[\ell(P^{J_{L'},L'})] - \IE^{\rho^{J_{l-1},l-1,\delta}}[\ell(P^{J_{L'},L'})] \right)- II_L
\;.
\end{eqnarray*}
Rearranging and using \eref{eq:PLtelescop} gives
\begin{eqnarray*}
\fl\IE^{\rho^{\delta}}[\ell(P)]=
&I_L+II_L+\sum_{l = 0}^L
\left( \IE^{\rho^{J_l ,l,\delta}}- \IE^{\rho^{J_{l-1},l-1,\delta}}\right)[\ell(P^{J_{L'},L'})]\\
&=I_L+II_L+
\sum_{l=0}^L\sum_{l' = 0}^{L'}
\left(\IE^{\rho^{J_l,l,\delta}}-\IE^{\rho^{J_{l-1},l-1,\delta}}\right)[\ell(P^{J
_{l'},l'}) - \ell(P^{J_{l'-1},l'-1})]
\;.
\end{eqnarray*}
The claimed expression follows from the definition of the MLMCMC method 
which computes the
MCMC sample path averages of the terms on the right hand side 
in \eref{eq:1zz}.
$\Box$ 

{
We next derive an error bound
by estimating}
the three terms in the error \eref{eq:MLMCMCerr} separately. 
Throughout we choose $J_l=2^{\lceil l/q\rceil }$.  
For the first term $I_L$, we obtain from Proposition \ref{prop:approxmu}, 
the bound
\beq\label{eq:ILbound}
|I_L| \leq C(\delta) (J_L^{-q} + 2^{-L}) \leq C(\delta) 2^{-L}
\;.
\eeq
Likewise we obtain, using Assumption \ref{assump:FEcomplexity} {in
addition to Proposition \ref{prop:approxmu}},
\beq\label{eq:IILbound}
|II_L| \leq C(\delta) \sum_{l=0}^L (J_\ell^{-q} + 2^{-l})2^{-L'(l)} 
       \leq C(\delta) \sum_{l=0}^L 2^{-(l+L'(l))}\;.
\eeq
%
We now estimate $III_L$. 
Since $J_l=2^{\lceil l/q\rceil}$ we have that 
\[
\sup_{u\in U}|\ell(P^{J_{l'},l'})-\ell(P^{J_{l'-1},l'-1})|\le C2^{-l'}\;. 
\]
Further
\beqas
\fl
&&\sup_{u\in U}|1-\exp(\Phi^{J_l,l}(u;\delta)-\Phi^{J_{l-1},l-1}(u;\delta))|
\\
&& \qquad
\le \sup_{u\in U}|\Phi^{J_l,l}(u;\delta)-\Phi^{J_{l-1},l-1}(u;\delta)|
(1+\exp(\Phi^{J_l,l}(u;\delta)-\Phi^{J_{l-1},l-1}(u;\delta)))
\\
&& \qquad
\le C2^{-l} 
\;.
\eeqas
As in Section \ref{sec:S},
for each discretization level $l\in \IN$, 
we introduce the Markov chains 
${{\mathcal C}_l=\{u^{(k)}\}_{k \in \IN_0}}\subset\mathbb{R}^{J_l}$ 
which are started in the restriction of $u^{(0)}$ 
to $\mathbb{R}^{J_l}$ and then
generated by \eref{eq:MC} with the 
$J_l$-term truncated parametric coefficient expansions \eref{eq:KtruncJ}
and forward problems \eref{eq:fe} which are Galerkin-discretized
at mesh level $l$, ie. with the acceptance probability 
$\alpha(u,v)$ in \eref{eq:alpha} 
replaced by
\be\label{eq:alphaJll}
\alpha^{J_l,l}(u,v)
=
1 \wedge \exp(\Phi^{J_l,l}(u;\delta)-\Phi^{J_l,l}(v;\delta))
\;,\quad (u,v)\in U\times U
\;.
\ee
For each discretization level $l=1,2,...$,
the chains ${\mathcal C}_l$ are pairwise independent. 
For every fixed discretization level $L$, 
we denote by 
${{\mathbf C}}_L =
 \{{\mathcal C}_{1},{\mathcal C}_{2},\ldots,{\mathcal C}_L\}$
the collection of Markov chains obtained from the different
discretizations,
and by ${{\mathbf P}}_L$
the product probability measure on the probability space
generated by the collection of these $L$ chains.
For each fixed discretization level $L$,
this measure describes the law of 
the collection of chains ${\mathbf C}_L$:
\[
{\mathbf P}_L
:=
\cP^{\rho,J_{1},1}\otimes \cP^{\rho,J_{2},{2}}\otimes\ldots\otimes\cP^{\rho,J_L,L}
\;.
\]
Let ${\mathbf E}_L$ be the expectation 
over all realizations of the collection ${\mathbf C}_L$ of chains 
    ${\mathcal C}_{l}$ with respect to the product measure ${\mathbf P}_L$. 
We have, 
by Assumption \ref{assump:FEcomplexity} 
and by Proposition \ref{t:mc},
\beqas
\fl &&{\mathbf E}_L\bigg[\bigg\{\IE^{\rho^{J_l,l,\delta}}
\left[(1-\exp(\Phi^{J_l,l}(u;\delta)-\Phi^{J_{l-1},l-1}(u;\delta)))](\ell(P^{J_{l'},l'})-\ell(P^{J_{l'-1},l'-1}))\right]
\\
\fl&&\qquad-E_{M_{ll'}}^{\rho^{J_l,l,\delta}}\left[(1-\exp(\Phi^{J_l,l}(u;\delta)-\Phi^{J_{l-1},l-1}(u;\delta)))](\ell(P^{J_{l'},l'})-\ell(P^{J_{l'-1},l'-1}))\right]\bigg\}^2\bigg]\\
\fl&&\qquad
\le C M_{ll'}^{-1}2^{-2(l+l')}
\;.
\eeqas
%
We also have, 
again using Assumption \ref{assump:FEcomplexity} and 
the last estimate in Proposition \ref{t:mc},
\beqas
\fl &&{\mathbf E}_L
\bigg[\bigg\{\IE^{\rho^{J_l,l,\delta}}[\exp(\Phi^{J_l,l}(u;\delta)-\Phi^{J_{l-1},l-1}(u;\delta))-1]\cdot\IE^{\rho^{J_{l-1},l-1,\delta}}[\ell(P^{J_{l'},l'})-\ell(P^{J_{l'-1},l'-1})]
\\
\fl &&\qquad-E_{M_{ll'}}^{\rho^{J_l,l,\delta}}[\exp(\Phi^{J_l,l}(u;\delta)-\Phi^{J_{l-1},l-1}(u;\delta))-1]\cdot E_{M_{ll'}}^{\rho^{J_{l-1},l-1,\delta}}[\ell(P^{J_{l'},l'})-\ell(P^{J_{l'-1},l'-1})]\bigg\}^2\bigg]
\\
\fl &&\le  
C\,{\mathbf E}_L
\bigg[\bigg\{\IE^{\rho^{J_l,l,\delta}}[\exp(\Phi^{J_l,l}(u;\delta)-\Phi^{J_{l-1},l-1}(u;\delta))-1]-E_{M_{ll'}}^{\rho^{J_l,l,\delta}}[\exp(\Phi^{J_l,l}(u;\delta)-\Phi^{J_{l-1},l-1}(u;\delta))-1]\bigg\}^2
\bigg]
\\
\fl &&\qquad
\cdot\sup_{u\in U}|\ell(P^{J_{l'},l'})-\ell(P^{J_{l'-1},l'-1})|^2+
\\
\fl && \qquad 
C\sup_{u\in U}|\exp(\Phi^{J_l,l}(u;\delta)-\Phi^{J_{l-1},l-1}(u;\delta))-1|^2
\\
\fl &&\qquad 
\cdot{\mathbf E}_L\bigg[\bigg\{\IE^{\rho^{J_{l-1},l-1,\delta}}
[\ell(P^{J_{l'},l'})-\ell(P^{J_{l'-1},l'-1})]-E_{M_{ll'}}^{\rho^{J_{l-1},l-1,\delta}}
[\ell(P^{J_{l'},l'})-\ell(P^{J_{l'-1},l'-1})]\bigg\}^2\bigg]
\\
\fl &&\le 
CM_{ll'}^{-1}2^{-2(l+l')}
\;.
\eeqas

Hence 
\beq\label{eq:IIILbound}
{\mathbf E}_L[|III_L|^2] 
\leq 
CL\sum_{l=0}^L L'(l)\sum_{l'=0}^{L'(l)} M_{ll'}^{-1} 2^{-2(l+l')} 
\;.
\eeq
To achieve a bound on the error \eref{eq:MLMCMCerr} which is 
$O(L^m2^{-L})$ for some $m\in \IR_+ $, 
we choose
\beq\label{eq:choicL'JM}
L'(l) := L - l , \quad \mbox{and}\quad M_{ll'} := 2^{2(L-(l+l'))}
\;. 
\eeq
We then have
\[
{\mathcal E}[|III_L|^2]
\le 
C L\sum_{l=0}^L(L-l)^22^{-2L}\le CL^42^{-2L}
\;.
\]
%
\begin{theorem}
\label{thm:MLMCMCerr}
For $d=2,3$, under the assumption that $P\in L^\infty(U,\rho;H^2(D)\bigcap V)$, 
and with the choices \eref{eq:choicL'JM} 
we have that
\beq\label{eq:MLMCMCerrest}
{\mathbf E}_L[| \IE^{\rho(\delta)}[P] - E^{MLMCMC}_L[P] |] 
\leq 
C(\delta)L^2 2^{-L}\;.
\eeq
The total number of degrees of freedom 
used in running the MLMCMC sampler is 
bounded by $O(L2^{2L})$ for $d=2$ and $O(2^{3L})$ for $d=3$.
Assuming the availability of a linear scaling multilevel 
iterative solver for the linear systems of equations arising
from Galerkin discretization, 
{the total number of floating point operations 
in the approximate evaluation of the
conditional expectation under the posterior 
(for one instance of data $\delta$)
is bounded} by $O(L^{d-1}2^{(d+1/q)L})$.
Denoting the total number of degrees of freedom 
which enter in running the chain on all discretization level 
by $N$, 
the error in \eref{eq:MLMCMCerrest} 
{is bounded by} $O((\log N)^{3/2}N^{-1/2})$ for $d=2$ and 
{by} $O((\log N)^{2}N^{-1/3})$ for $d=3$.
The total number of floating point operations {to realize the MLMCMC} 
is bounded by $O((\log N)^{-1/(2q)}N^{1+1/(2q)})$ for $d=2$ and  {by}
$O((\log N)^{2}N^{1+1/(3q)})$ for $d=3$. 
\end{theorem} 
{\it Proof}\ \ 
At each step we solve the truncated forward equation 
\eref{eq:probJ} 
for the truncation levels $J_l$ and $J_{l'}$ in expansion 
\eref{eq:KtruncJ} 
and the Finite Element discretization levels $l$ and $l'$, respectively. 
Assuming approximate solution of the discretized problems
with termination at the discretization error $O(2^{-l})$, resp. $O(2^{-l'})$,
and a linear scaling iterative solver with the choices \eref{eq:choicL'JM}
the total work is bounded by an absolute multiple of 
the total number of degrees of freedom,
which in turn is bounded, for $d\geq 2$,
by
\beqas
\sum_{l=0}^L\sum_{l'=0}^{L'(l)}M_{ll'}(2^{dl}+2^{dl'})
=
2^{2L}\sum_{l=0}^L
\sum_{l'=0}^{L'(l)}\left(2^{(d-2)l}\cdot 2^{-2l'}+2^{-2l}\cdot 2^{(d-2)l'}\right)
\\
\lesssim 2^{2L}\left(\sum_{l=0}^L2^{(d-2)l}+\sum_{l=0}^L2^{-2l}\sum_{l'=0}^{L-l}2^{(d-2)l'}\right).
\eeqas
Therefore for $d=2$, the number of degrees of freedom is bounded by
\[
\lesssim 2^{2L}\left(L+\sum_{l=0}^L2^{-2l}(L-l)\right)\lesssim L2^{2L}.
\]
For $d=3$, it is bounded by
\[ 
\lesssim 2^{2L}\left(2^{L}+\sum_{l=0}^L2^{-2l}2^{L-l}\right)
\lesssim 2^{2L}\left(2^{L}+\sum_{l=0}^L2^{L}2^{-3l}\right)
\lesssim 2^{{3}L}
\;.
\]
From the proof of Proposition \ref{prop:fecost}, 
we infer that the number of floating point operations at each step is 
not larger than $O(l^{d-1}2^{dl}J_l+{l'}^{d-1}2^{dl'}J_{l'})$ so that, 
using again \eref{eq:choicL'JM},
the total number of floating point operations required
to evaluate the MLMCMC estimator is bounded by 
\beqas
\lesssim\sum_{l=0}^L\sum_{l'=0}^{L'(l)}M_{ll'}(l^{d-1}2^{dl}J_l+{l'}^{d-1}2^{dl'}J_{l'})
\\
\lesssim 2^{2L}\sum_{l=0}^L\sum_{l'=0}^{L'(l)}
(l^{d-1}2^{(d-2+1/q)l}2^{-2l'}+l'^{d-1}2^{(d-2+1/q)l'}2^{-2l})
\\
\lesssim L^{d-1}2^{2L}\left(\sum_{l=0}^L2^{(d-2+1/q)l}+\sum_{l=0}^L2^{(d-2+1/q)(L-l)}2^{-2l}\right) 
\\
\lesssim L^{d-1}2^{(d+1/q)L}
\;.
\eeqas
\hfill$\Box$
\begin{remark} 
\label{remk:ON-MLsolver}
In the preceding work analysis we assumed
the availability of 
linear scaling, iterative solvers 
for the approximate solution of 
the Finite Element discretizations \eref{eq:fe} 
of the parametric, elliptic problems \eref{eq:probJ}. 
These include Richardson or conjugate gradient iteration
with diagonal preconditioning using 
Riesz bases for the Finite Element spaces $V^l$ 
as presented in Appendix D. 
However, if such bases are not available, 
our results will also apply to standard Finite Element discretizations,
provided that linear scaling multilevel solvers, such as multigrid, 
are used to solve \eref{eq:fe}, for each increment of the Markov chain.
Our work versus accuracy analysis will also apply in these cases, 
with identical conclusions.
\end{remark}

\section{Conclusions}
\label{sec:Concl}
Rewriting the convergence bounds 
in Theorems \ref{t:fpe}, \ref{thm:improvedMCMC} 
and \ref{thm:MLMCMCerr} in terms of the error, 
to achieve an error $\ep$ in the mean square with 
respect to the probability measure 
that generates the randomness of the Markov Chain 
ie., ${\cal P}^{\rho,J,l}$, ${\cal P}^{\rho,N,{\cal L}}$ and ${\mathbf P}_L$, 
respectively 
we bound the complexity as follows (ignoring the logarithmic factors for conciseness):
\begin{itemize}
\item for the Plain MCMC procedure: $O(\ep^{-d-2})$  degrees of freedom and 
       $O(\ep^{-d-2-1/q})$ floating  point operations;

\item for the gpc-MCMC procedure: $O(\ep^{-1/\tau})$ degrees of freedom 
      and $O(\ep^{-\max(\alpha/\tau,2+1/\sigma)})$ floating point operations;

\item and for the MLMCMC procedure:   
      the essential optimal complexity of 
      $O(\ep^{-d})$ degrees of freedom and $O(\ep^{-d-1/q})$ floating point operations.
\end{itemize}
Therefore, 
the complexity of the gpc-MCMC is superior to that of the Plain MCMC 
procedure when in Assumption \ref{ass:gpcplxity} $\tau$ is close to
$1/d$ and $\alpha$ is close to $1$, which holds in many cases 
(see, e.g. \cite{ScCJG} and the references therein). 
We also have shown that the asymptotic accuracy vs. work 
of the MLMCMC is always superior to that of the Plain MCMC. 
We note that the complexity of the MLMCMC is of the same order 
as that for solving one single realization of equation \eref{eq:prob} 
with the same level of truncation for the coefficient, which indicates
some form of optimality of the MLMCMC method proposed here.

We have considered the MLMCMC estimation of functionals $\ell(\cdot)$ of
the solution; in doing so, we used only minimal regularity $\ell(\cdot)\in V^*$
of these functionals, leading to the convergence rates \eref{eq:1}. 
Higher regularity $\ell(\cdot)\in L^2(D)$ 
will entail in \eref{eq:1} the Finite Element convergence rate $2^{-2L}$, 
via a classical Aubin-Nitsche duality argument.
Likewise, we only used the lowest order Finite Element methods. 
In the present analysis, we focused on first order Finite Element discretizations 
as stronger regularity assumptions on $f$, on $\ell(\cdot)$
and on the $\psi_j$ will {\em not} imply corresponding higher 
rates of convergence for gpc-MCMC and ML-MCMC, 
due to the (maximal) order $1/2$ afforded by the MC method.

A further aspect pertaining to the overall complexity is the
following: in \cite{SchwStua12} an
entirely deterministic approach to the solution
of the Bayesian inverse problem is presented,
and in \cite{Schil813} corresponding numerical 
experiments are presented. 
These methods offer convergence rates in terms of 
the number $M$ of samples
which are superior to the rate $1/2$ 
of MC methods such as those analyzed here.
{Nonetheless, MCMC methods will remain popular because
of their data-adaptive nature; the {present results
indicate how the use of}
gpc and multilevel ideas {may be used to attain
significant speedups of MCMC methods}.
}
%

We also observe that we have concentrated on a very special
MCMC method, namely the independence sampler. This has been dictated
by the need to use MCMC methods which scale-independently of
dimension and {\em for which there is a complete 
analysis of the convergence resulting chain} demonstrating this fact.
Whilst there are now several dimension-independent MCMC methods
\cite{cotter2012mcmc,vollmer2013} the independence sampler
is the only one for which the required analysis of the convergence
properties is sufficiently developed for our theory; we
anticipate further theoretical developments for different MCMC
methods, and different inverse problems, in the future. However
we do note that the independence sampler will work well
when the negative log likelihood $\Phi$ does not vary significantly,
although it will be inefficient in general. For problems where
the prior is Gaussian or log-normal Gaussian,
appropriate MCMC methods may be found in \cite{cotter2012mcmc}. 
However for these more general methods the analysis of the Markov chain
based on the methods of \cite{MT93}, {as we have used
for the independence sampler here}, is not appropriate
and more sophisticated arguments are required, 
as presented in \cite{hairer2013}. 

\section*{Appendix A: Lipschitz Dependence of the Forward Map on Parameters}

{\it Proof of Proposition \ref{prop:measurable}\ \ }
The existence of a solution $P$ to \eref{eq:prob}, and the
bound \eref{eq:stndLaxMilg}, follows from standard application
of Lax-Milgram theory. We have, for every $\phi\in V$, $u,u'\in U$,
\begin{eqnarray}
\int_DK(x,u)(\nabla P(x,u)-\nabla P(x,u'))\cdot\nabla\phi(x){dx}
\nonumber
\\
=\int_D(K(x,u')-K(x,u))\nabla P(x,u')\cdot\nabla\phi(x)dx \;. 
\label{eq:diff}
\end{eqnarray}
Again using {\eref{eq:KminKmax}}, i.e. 
that $K(x,u)$ is bounded below uniformly 
with respect to $(x,u) \in D \times U$, 
it follows that there exists $C>0$ 
such that for all $u\in U$
\be
\|P(\cdot,u)-P(\cdot,u')\|_{V} 
\le 
C\| P(\cdot,u')\|_V\|K(\cdot,u')-K(\cdot,u)\|_{L^\infty(D)}
\;.
\label{eq:Lipschitzofp1}
\ee
Due to \eref{eq:stndLaxMilg},
it follows from \eref{eq:Lipschitzofp1}
that there exists a constant $C>0$
such that
\be
\forall u\in U:\quad 
\|P(\cdot,u)-P(\cdot,u')\|_{V} 
\le 
C\|K(\cdot,u')-K(\cdot,u)\|_{L^\infty(D)}
\;.
\label{eq:Lipschitzofp2}
\ee
From \eref{eq:a} and Assumption \ref{assump:sumpsi}(i) 
it follows with $C>0$ as in \eref{eq:Lipschitzofp2}
that
\begin{eqnarray*}
\|P(\cdot,u)-P(\cdot,u')\|_{V} 
&\le& C \sum_{\geq 1} |u_j-u_j'|\|\psi_j\|_{L^\infty(D)}\;
\\
&\le& C \|u-u'\|_{\ell^{\infty}(\IN)}\sum_{j\geq 1} \|\psi_j\|_{L^\infty(D)}\;
\\
&\le& C
{\kappa\over 1+\kappa}\bar K_{\min}
\|u-u'\|_{\ell^{\infty}(\IN)}
\;.
\end{eqnarray*}
This establishes the desired \Lipsh continuity,

\section*{Appendix B: Bayesian inverse problems in measure spaces}
\label{sec:B}
On a {measurable} space $(U,\Theta)$ 
where $\Theta$ is a $\sigma$-algebra 
consider a measurable map $\cG:\,U\to (\spa,\cB^k)$. 
The data $\delta$ is assumed to be an 
observation of $\cG$ subject to 
an unbiased observation noise $\vartheta$:
\[
\delta = \cG(u) + \vartheta.
\]
We assume that $\vartheta$ is a centred Gaussian 
with law $N(0,\Sigma)$. 
Let $\rho$ be a prior
probability measure on $(U,\Theta)$. 
Our purpose is to 
determine the conditional probability $\IP(u|\delta)$ 
on $(U,\Theta)$. 
The following result holds.
\begin{proposition}\label{prop:RN}
Assume that $\cG:U\to\spa$ is measurable. 
The posterior measure $\rho^\delta(du) = \IP(du|\delta)$ 
given data $\delta$
is absolutely continuous with respect to the prior measure 
$\rho(du)$ and has the Radon-Nikodym derivative
\eref{eq:posterior} with $\Phi$ given by \eref{eq:phi}.
\end{proposition}
This result is established in Cotter et al.\cite{Cotteretal} and Stuart \cite{AndrewActa}. 
Though the setting in  \cite{Cotteretal} and \cite{AndrewActa} is in a Banach space $X$, 
the proofs of Theorem 2.1 in \cite{Cotteretal} and Theorem 6.31 of \cite{AndrewActa} 
hold for any measurable spaces as long as the mapping $\cG$ is measurable. 

To study the well-posedness of the posterior measures, that is
continuity with respect to changes in the observed data,
we use the Hellinger distance, as in Cotter et al. 
\cite{Cotteretal}, 
which is defined as
\be
d_{\rm{Hell}}(\mu,\mu')
=
\left(\frac{1}{2}\int_U\left(\sqrt{{d\mu\over d\rho}}-\sqrt{d\mu'\over d\rho}\right)^2d\rho\right)^{1/2}
\label{eq:Hell}
\ee
for any two measures $\mu$ and $\mu'$ on $U$ 
which are absolutely continuous with respect to a common
reference measure $\rho$\footnote{
Note that {\em any} such common
reference measure will deliver the same value for the Hellinger
distance}.
In \cite{Cotteretal},
it is proved that when $U$ 
is a  Banach space, if the prior measure $\rho$ is Gaussian, 
and under the conditions that $\Phi$ grows polynomially with respect to $u$, 
and is locally \Lipsh with respect to $u$ fixing $y$ and with respect to $y$ fixing $u$,  
in the second case with a \Lipsh constant which also grows polynomially in $u$,
then the posterior measure given the data $\delta$, 
i.e.  $\rho^\delta$,
is locally \Lipsh in the Hellinger distance $d_{\rm {Hell}}$:
\[
d_{\rm {Hell}}(\rho^\delta, \rho^{\delta'})\le c |\delta - \delta'|\;,
\] 
where (recall) 
$|\cdot|$ denotes the Euclidean distance in $\spa$. 
The Fernique theorem plays an essential role in the proofs,
exploiting the fact that polynomially growing
functions are integrable with respect to Gaussians.
In this section, we extend this result to measurable spaces under more 
general conditions than those in Assumption 2.4 of Cotter et al. \cite{Cotteretal}; 
in particular we do not assume a Gaussian prior.
{
The following assumption concerning
the local boundedness, 
and local \Lipsh dependence
of $\Phi$ on $\delta$, will be crucial in our argument.
Its validity for the model problem \eref{eq:prob}, 
with \eref{eq:a} and Assumption \ref{assump:sumpsi}
will be verified in the ensuing proof of Proposition \ref{prop:new}.
}
\begin{assumption} \label{assump:Phi}
Let $\rho$ be a probability measure on the measure
space $(U,\Theta)$.
The Bayesian potential function 
$\Phi:U\times \spa\to {\mathbb R}$ 
satisfies:
\begin{itemize}
\item[(i)] 
(local boundedness)
for each $r>0$, there is a constant $\eM(r)>0$ and 
a set $U(r)\subset U$ 
of positive $\rho$ measure such that for all $u\in U(r)$ and 
for all $\delta $ such that $|\delta|_{\Sigma} \le r$
\be\label{eq:boundPhi}
0 \le \Phi(u;\delta) \le \eM(r)\;.
\ee
\item[(ii)] 
(local Lipschitz continuity of Bayesian Potential $\Phi$ on data $\delta$)
there is a mapping 
$G:{\mathbb R}\times U\mapsto {\mathbb R}$ 
such that for each $r>0$, $G(r,\cdot)\in L^2(U,\rho)$; 
and 
for every $|\delta|_\Sigma , |\delta'|_\Sigma \le r$ 
it holds that
\[
|\Phi(u;\delta) - \Phi(u;\delta')| \le G(r,u)|\delta - \delta'|_\Sigma
\;.
\]
\end{itemize}
\end{assumption}
Under Assumption \ref{assump:Phi}, 
the definition \eref{eq:posterior} of the 
posterior measure $\rho^\delta$ is meaningful
as we now demonstrate. 

\begin{proposition}\label{prop:wellposedness} 
Under Assumption \ref{assump:Phi}, 
the measure $\rho^\delta$ 
depends {locally} \Lipsh continuously 
on the data $\delta$ with respect to the Hellinger metric: 
for each positive constant $r$ there 
is a positive constant $C(r)$ such that 
if $|\delta|_\Sigma,\,|\delta'|_\Sigma\le r$, 
then
\[
d_{\rm{Hell}}(\rho^\delta,\rho^{\delta'}) 
\le C(r)|\delta - \delta'|_\Sigma 
\;.
\]
\end{proposition}
 {\it Proof} \ \ 
Throughout this proof $K(r)$ denotes a constant
depending on $r$, possibly changing from instance to
instance.
The normalization constant in \eref{eq:posterior} is 
\be
Z(\delta) = \int_U \exp(-\Phi(u;\delta))d\rho(u)\;.
\label{eq:zy}
\ee
We first show that for each $r>0$, there is a positive
constant $K(r)$ such that $Z(\delta) \ge K(r)$ 
when $|\delta |_\Sigma \le r$. 
From (\ref{eq:zy}) and 
Assumption \ref{assump:Phi}(i) it follows 
that when $|\delta |_\Sigma \le r$, there holds
\beq  \label{eq:Zboundlow}
Z(\delta)\ge \rho(U(r))\exp(-\eM(r)) > 0\;.
\eeq
Using the inequality 
$|\exp(-x)-\exp(-y)|\le |x-y|$ which holds for $x,y \ge0$ we find
\[
|Z(\delta)-Z(\delta')| \le \int_U|\Phi(u;\delta) - \Phi(u;\delta')|d\rho(u)
\;.
\]
From Assumption \ref{assump:Phi}(ii),
\beqas
|\Phi(u;\delta ) - \Phi(u;\delta')| \le G(r,u)|\delta - \delta'|_\Sigma 
\;.
\eeqas
As $G(r,u)$ is $\rho$-integrable, 
there is $K(r)$ such that 
\[
|Z(\delta) - Z(\delta')| \le K(r)|\delta - \delta'|_\Sigma 
\;.
\]
The Hellinger distance satisfies
\begin{eqnarray*}
\fl 2d_{\rm{Hell}}(\rho^\delta ,\rho^{\delta'})^2 
&=& \int_U\Bigl(Z(\delta)^{-1/2}\exp(-{1\over 2}\Phi(u;\delta))\;-
Z(\delta')^{-1/2}\exp-{1\over 2}\Phi(u;\delta')\Bigr)^2d\rho(u)
\\
&\le& I_1+I_2,
\end{eqnarray*}
where 
\[
I_1={2\over Z(\delta)} 
\int_U\Bigl(\exp(-{1\over 2}\Phi(u;\delta))-\exp(-{1\over 2}\Phi(u;\delta'))\Bigr)^2d\rho(u),
\]
and
\[
I_2 = 2|Z(\delta)^{-1/2}-Z(\delta')^{-1/2}|^2 \int_U\exp(-\Phi(u;\delta'))d\rho(u).
\]
Using again $|\exp(-x)-\exp(-y)|\le |x-y|$, 
we have, for constant $K(r)>0$,
\beqas
I_1\le K(r)\int_U|\Phi(u;\delta) - \Phi(u;\delta')|^2d\rho(u)
\\
\le K(r)\int_U(G(r,u))^2 d\rho(u)|\delta - \delta'|_\Sigma^2 
\le 
K(r)|\delta - \delta'|_\Sigma^2
\;.
\eeqas 
Furthermore,
\[
|Z(\delta)^{-1/2} - Z(\delta')^{-1/2}|^
2\le
K(r)|Z(\delta) - Z(\delta')|^2 \le K(r)|\delta - \delta'|_\Sigma^2\;.
\]
The conclusion follows.
\hfill$\Box$

\begin{proposition}
\label{prop:new}
For the elliptic PDE \eref{eq:prob}, the function $\cG$ defined
by \eref{eq:defG} and viewed as map from 
$U$ to ${\mathbb R}^k$ is \Lipsh, 
if $U$ is endowed with the topology of $\ell^{\infty}(\IN)$.
Moreover, Assumption \ref{assump:Phi}
holds with $U(r)=U$ and $G(r,u)=G(r)$ independent of $u$. 
\end{proposition}
 {\it Proof} \ \ 
We have
\[
\forall u,u'\in U:\quad
|\cG(u)-\cG(u')|
\le 
C\max_{i}\{\|\cO_i\|_{V^*}\}\|P(\cdot,u)-P(\cdot,u')\|_{V}\;.
\]
From \eref{eq:Lipschitzofp2} there exists a constant $c>0$
such that 
\[
\forall u,u'\in U:\quad
|\cG(u)-\cG(u')|\le 
C\|K(\cdot,u) - K(\cdot,u')\|_{L^\infty(D)}
\;.
\]
From Proposition \ref{prop:measurable}, 
we deduce that $\cG$ as map from $U \subset \ell^{\infty}(\IN)$ 
to ${\mathbb R}^k$ is \Lipsh. 

We now verify Assumption \ref{assump:Phi}. 
For the function $\cG(u)$ we have from \eref{eq:defG}
for every $u\in U$ the bound
\[
|\cG(u)| \le \max_i\{\|\cO_i\|_{V^*}\}\|P(\cdot,u)\|_{V}
\;.
\]
From \eref{eq:stndLaxMilg}, $\sup \{ |\cG(u)| : u\in U \} < \infty$.
We note that for given data $\delta$, there holds
\[
\forall u\in U:\quad 
|\Phi(u;\delta)| \le {1\over 2}(|\delta|_\Sigma + |\cG(u)|_\Sigma)^2
\]
and hence, since $\sup_{u\in U} |\cG(u)|$ is finite,
the set $U(r)$ in Assumption \ref{assump:Phi}(i) 
can be chosen as $U$ for all $r$.
We have, for every $u\in U$, 
\beqas
|\Phi(u;\delta) - \Phi(u;\delta')|
\le 
{1\over 2}|\langle\Sigma^{-1/2}(\delta+\delta'-2\cG(u)),\Sigma^{-1/2}(\delta - \delta')\rangle|
\\
\le 
{1\over 2}
(|\delta|_\Sigma + |\delta'|_\Sigma + 2|\cG(u)|_\Sigma)|\delta-\delta'|_\Sigma 
\;.
\eeqas
Choosing the function $G(r,u)$ in Assumption \ref{assump:Phi}(ii) 
as
\[
G(r,u)={1\over 2}(2r+C),
\]
for a sufficiently large constant $C>0$ 
(depending only on bounds of
the forward map and the observation functionals, but independent of
the data $\delta$ and of $u$), 
we have shown that Assumption \ref{assump:Phi}(ii) holds in
the desired form. 
With Proposition \ref{prop:wellposedness} follows the assertion.
\hfill$\Box$

{\it Proof of Proposition \ref{p:wp}}\ \ 
From Proposition \ref{prop:new}, 
we deduce that $\cG$ as map from $U \subset \ell^{\infty}(\IN)$ 
to ${\mathbb R}^k$ is \Lipsh and, hence,
$\rho$-measurable. 
We then apply Proposition \ref{prop:RN} to
deduce the existence of $\rho^\delta$ and the formula for
its Radon-Nikodym derivative with respect to  $\rho$. 
Proposition \ref{prop:wellposedness} gives the desired
Lipschitz continuity, since Proposition \ref{prop:new} establishes
Assumption \ref{assump:Phi}.
\hfill$\Box$
\section*{Appendix C: Convergence Properties of the Independence Sampler}
{\it Proof of Proposition \ref{t:mc}\ \ }
We claim that \eref{eq:MC} 
defines a Markov chain $\{u^{(k)}\}_{k=0}^{\infty}$ which is reversible
with respect to $\rho^\delta$. 
To see this let
$\nu(du,dv)$ denote the product measure 
$\rho^\delta(du)\otimes \rho(dv)$ 
and 
$\nu^{\dagger}(du,dv)=\nu(dv,du)$.
Note that $\nu$ describes the probability distribution
of the pair $(u^{(k)},v^{(k)})$ on $U \times U$ when
$u^{(k)}$ is drawn from the posterior distribution
$\rho^\delta$, and $\nu^{\dagger}$
designates the same measure 
with the roles of $u$ and $v$ reversed.
These two measures are equivalent as measures
if $\rho^\delta$ and $\rho$ are equivalent, which we establish
below; it then follows that
\begin{equation}
\label{eq:RND2}
\frac{d\nu^{\dagger}}{d\nu}(u,v) 
= 
\exp\bigl(\Phi(u;\delta ) - \Phi(v;\delta) \bigr)\;,
\quad (u,v)\in U\times U 
\;.
\end{equation} 
From Proposition 1 and Theorem 2 in \cite{T98} 
we deduce that \eref{eq:MC}
is a Metropolis-Hastings Markov chain which is
$\rho^\delta$ reversible, since $\alpha(u,v)$ given by
\eref{eq:alpha} is equal to 
$\min\{1, \frac{d\nu^{\dagger}}{d\nu}(u,v)\}$.

Equivalence of $\rho^\delta$ and $\rho$ follows since 
the negative of the log-density is bounded from above and
below, uniformly on $U$, because Proposition \ref{prop:new}
shows that \eref{eq:boundPhi} holds with $U(r)=U$.
Using \eref{eq:boundPhi} and \eref{eq:alpha} 
it follows that the proposed random draw from $\rho$ has probability
greater than $\exp\bigl(-\eM(r)\bigr)$ of being accepted. 
Thus
$$
p(u,A) \ge \exp\bigl(-\eM(r)\bigr)\rho(A) \quad \forall 
u \in U\;.
$$ 
The first result follows from \cite{MT93}, Theorem 16.2.4 
with $X=U.$ 
The second result follows from \cite{MT93}, Theorem 17.0.1. 
To see that the constant $c$ in \eref{eq:ex}
can be bounded only in terms of
$\eM(r)$ and $\sup_{u\in U}|g(u)|$,
we note that it is given by 
\be\label{eq:defc2}
c^2 
=
\cE^{\rho^\delta}|\bar g(u^{(0)})|^2
+
2\sum_{n=1}^\infty\cE^{\rho^\delta}[\bar g(u^{(0)})\bar g(u^{(n)})]
\ee
where to ease notation we introduced the function $\bar g$ 
as $\bar g = g-\IE^{\rho^\delta}(g)$
(and we do not tag the dependence of $\bar g$ on the data
     $\delta$ in the remainder of this proof).
Note that $\sup_{u\in U} |\bar g (u) | $
      is bounded uniformly w.r. to the data $\delta$. 
The equation \eref{eq:defc2} is commonly known as the 
integrated autocorrelation time. Now
\beqas
2\sum_{n=0}^\infty{\cE^{\rho^\delta}}[\bar g(u^{(0)})\bar g(u^{(n)})]
&\le&2\sup_u|\bar g(u)|\IE^{\rho^\delta}
\sum_{n=0}^\infty
|{\cE_{u(0)}}[{\bar g}(u^{(n)})]|
\\
& = &
2\sup_u|\bar g(u)|\IE^{\rho^\delta}
\sum_{n=0}^\infty
|{\cE_{u(0)}}[g(u^{(n)})]-\IE^{\rho^\delta}[g]|
\\
&\le &
4\sup_u|\bar g(u)|^2
\sum_{n=0}^\infty\Bigl(1-\exp\bigl(-\eM(r)\bigr)\Bigr)^n
\;.
\eeqas 
For the mean square approximation, using
the stationarity of the Markov chain 
conditioned to start in $U \ni u^{(0)} \sim \rho^\delta$,
we have
\beqas
\fl{1\over M}{\cE^{\rho^\delta}}\Big[\Big|\sum_{k=1}^M\bar g(u^{(k)})\Big|^2\Big]
&=&
\IE^{\rho^{\delta}}[\bar g(u^{(0)})^2]+2{1\over M}
\sum_{k=1}^M\sum_{j=k+1}^M\cE^{\rho^\delta}[\bar g(u^{(k)})\bar g(u^{j})]
\\
&=& \ds
\IE^{\rho^{\delta}}[\bar g(u^{(0)})^2]
+
2{1\over M}\sum_{k=0}^{M-1}\sum_{j=1}^{M-k}\cE^{\rho^ \delta}[\bar g(u^{(0)})\bar g(u^{(j)})]
\\
&=& \ds
\IE^{\rho^\delta}[\bar g(u^{(0)})^2]
+
2{1\over M}
\sum_{k=0}^{M-1}\sum_{j=1}^{M-k}\IE^{\rho^{\delta}}[\bar g(u^{(0)})\cE_{u^{(0)}}[\bar g(u^{(j)})]]
\\ 
&\le& \ds
\IE^{\rho^{\delta}}[\bar g(u^{(0)})^2]
\\ 
& & \ds
+2{1\over M}\sum_{k=0}^{M-1}{\sup_{u}|\bar g(u)|}
\sum_{j=1}^{M-k}\IE^{\rho^\delta}[|\cE_{u^{(0)}}[g(u^{(j)})]-\IE^{\rho^\delta}[g]|]
\\
&\le& \ds
\IE^{\rho^{\delta}}[\bar g(u^{(0)})^2]
+
4{1\over M}\sum_{k=0}^{M-1}\sup_{u}|\bar g(u)|^2 \sum_{j=1}^{M-k} \Bigl(1-\exp\bigl(-\eM(r)\bigr)\Bigr)^j\\ 
&\le& \ds
\IE^{\rho^{\delta}}[\bar g(u^{(0)})^2]
+
4\sup_{u}|\bar g(u)|^2 \sum_{j=1}^{\infty} {\Bigl(1-\exp\bigl(-\eM(r)\bigr)\Bigr)^j}, 
\eeqas
which is clearly bounded uniformly with respect to $M$. 
Thus we have shown that there exists $C>0$ such that for all $M$
\[
{\cE^{\rho^\delta}}\Big[\Bigl|\frac{1}{M}\sum_{k=1}^M\bar g(u^{(k)})\Big|^2\Big]\le \frac{C}{M}
\;.
\]
It remains to show that the expectation 
{$\cE^{\rho^\delta}$}
with respect 
to the unknown posterior $\rho^{\delta}$ 
can be replaced by an expectation 
with respect to the prior measure $\rho$.

To this end we note that 
\beqas
\cE^{\rho}\Big[\Big|\sum_{k=1}^M\bar g(u^{(k)})\Big|^2\Big]
&=&\ds
\int_U\cE_{u^{(0)}}\Big[\Big|\sum_{k=1}^M\bar g(u^{(k)})\Big|^2\Big]d\rho(u^{(0)})\\
&=&\ds
\int_U
\cE_{u^{(0)}}\Big[\Big|\sum_{k=1}^M\bar g(u^{(k)})\Big|^2\Big]{d\rho\over d\rho^\delta}(u^{(0)})d\rho^\delta(u^{(0)})
\\
&\le&\ds
\cE^{\rho^\delta}\Big[\Big|\sum_{k=1}^M\bar g(u^{(k)})\Big|^2\Big] 
Z(\delta)\sup_{u\in U}\Big[\exp(\Phi(u;\delta))\Big]
\;.
\eeqas
As $Z(\delta) \le 1$ and $\Phi(\cdot;\delta)$ 
is assumed to be bounded uniformly,
we deduce that
\[ 
\cE^{\rho}\Big[\Big|\frac{1}{M}\sum_{k=1}^M\bar g(u^{(k)})\Big|^2\Big]\le \frac{C}{M},
\]
for a constant $C$ independent of $M$. 
The conclusion then follows.  \hfill$\Box$

\section*{Appendix D: Finite Element Methods}
%
In this Appendix we prove that Assumption \ref{assump:FEcomplexity} holds
if we employ a Riesz basis for the Finite Element space and 
when the domain and the equation's coefficients possess sufficient regularity.

We assume in the following that
{\em 
the union of all Finite Element basis functions $w^{l}_j$ 
of the spaces $V^l={\rm span}\{ w^l_j : j=1,...,N_l\}$, 
$l=0,1,2,...$,
constitutes a  Riesz basis in $V$}.
We remark that such bases are available 
in two and three dimensional polyhedral domains
(see, e.g., \cite{NguyenStevenson})
(the following assumption of availability of 
$V$-stable Riesz bases is made for convenience,
and may also be replaced by 
the availability of a linear complexity, optimal
preconditioning, such as a BPX preconditioner).

\begin{assumption}\label{assump:wavelet} 
(Riesz Basis Property in $V$)
For each $l\in\IN_0$ there exists a set of indices 
$I^l\subset\IN^d$ of cardinality $N_l = O(2^{ld})$ 
and a family of basis functions 
$w^{l}_k\in H^1_0(D)$ 
indexed by a multi-index $k\in I^l$ 
such that 
$V^l={\rm span}\{w^{l}_k \ : k\in I^l\}$, 
and there exist constants $c_1$ and $c_2$ 
which are independent of the discretization
level $l$ 
such that if $w\in V^l$ is written as
$w = \sum_{k\in I^l} c^{l}_k w^l_k \in V^l$, 
then
\[
c_1\sum_{k\in I^l}|c^l_k|^2
\le
\|w\|_{V}^2
\le 
c_2\sum_{k\in I^l}|c^{l}_k |^2
\;.
\]
\end{assumption}
Multiscale Finite Element bases entail, in general, 
larger supports than the standard, single scale
basis functions which are commonly used in the Finite 
Element Method, which implies that the stiffness matrices
in these bases have additional nonzero entries, as 
compared to $O({\rm dim} V^l) = O(2^{dl})$ 
many nonzero entries
of the stiffness matrices that result when 
one-scale bases, such as the hat functions,
are used.

To bound the number of nonzero entries, 
we shall work under 
\begin{assumption}\label{assump:support}
(Support overlap)
For all $l\in \IN_0$ and for every $k\in I^l$,
for every $l'\in \IN_0$ the support intersection
$\supp(w^{l}_k)\cap\supp(w^{l'}_{k'})$ 
has positive measure for at most 
$O(\max(1,2^{l'-l}))$ values of $k'$.
\end{assumption}
Assumption \ref{assump:support} implies that the number of 
non-zero entries in the stiffness matrix of the approximating problem \eref{eq:fe} is 
at most $O(l^{d-1}2^{dl})$. 
To prove the error bound \eref{eq:fe2},
we require the regularity $P(\cdot,u)$ belonging and 
being bounded in $H^2(D)\bigcap H^1_0(D)$ 
{uniformly with respect to $u\in U$.}
Assumption \ref{assump:sumpsi}(ii) implies
the following regularity results. 
\begin{proposition}\label{prop:H2}
If $D$ is convex and $f \in L^2(D)$, 
and if Assumptions \ref{assump:sumpsi}(i),(ii) hold,
then, for every $u\in U$,
the solution $P^J(\cdot,u)$ of \eref{eq:probJ} 
belongs to 
the space 
$W := H^2(D)\cap H^1_0(D)$ 
and there exists a positive constant $C>0$
such that 
$$
\sup_{J\in \IN} \sup_{u\in U} \| P^J(\cdot,u) \|_W \leq C \| f \|_{L^2(D)} \;.
$$
\end{proposition}
{\it Proof}\ \ 
By \eref{eq:KminKmax},  $K^J(x,u) \geq K_{min} > 0$ and
we may rewrite the PDE in \eref{eq:probJ} 
as
\[
-\Delta P^J(x,u) = {1\over K^J(x,u)}(f(x)+\nabla K^J(x,u)\cdot\nabla P^J(x,u)).
\]
By our assumptions, 
the right hand side is uniformly bounded with respect to $J$ and
$u\in U$ in the space $L^2(D)$.
As the domain $D$ is convex, we deduce that $P^J$ are uniformly bounded
with respect to $J$ and $u\in U$ in the space $W$: it holds
$$
\begin{array}{rcl}
\sup_{u\in U} \|\Delta P^J(\cdot,u) \|_{L^2(D)}
&\leq & \ds 
\frac{1}{K_{\min}} 
\sup_{u\in U} \sup_{J\geq 1}
\left[
\| f \|_{L^2(D)}
+
\| K^J (\cdot,u) \|_{W^{1,\infty}(D)}
\| P^J (\cdot,u) \|_V
\right]
\\
&\leq& \ds C < \infty \;
\end{array}
$$
due to the summability of the $W^{1,\infty}(D)$-norms
of the $\psi_j$ implied by Assumption \ref{assump:sumpsi}(ii).
The desired, uniform (w.r. to $J$ and $u$) 
bound in the $W$ norm then follows from the 
$L^2$ bound on $\Delta P^J(\cdot,u)$ 
and \eref{eq:stndLaxMilg} and the convexity of the domain $D$.
\hfill$\Box$

\section*{Appendix E: Generalized Polynomial Chaos Methods}
We justify Assumption \ref{ass:gpcplxity} in this section. 
For the ensuing analysis, we shall 
impose the following assumption on the 
summability of the gpc expansion of $P$:
\begin{assumption}\label{ass:gpclp}
There exists a constant $0<p<1$ 
such that the coefficients $P_\nu$ of the gpc expansion of $P$ 
satisfy $(\|P_\nu\|_V)_\nu\in\ell^p({\mathcal F})$.
\end{assumption}
This assumption is valid under the provision of suitable
decay of the coefficient functions $\psi_j$ such as 
Assumption \ref{assump:sumpsi}(ii). 
We refer to \cite{CDS1,CDS2} for details.
By a classical argument (``Stechkin's Lemma''),
this implies the following, 
so-called ``best $N$-term approximation property''.
\begin{proposition}\label{prop:Stechkin}
Under Assumption \ref{ass:gpclp},
there exists a nondecreasing sequence
$\{ \Lambda_N \}_{N\in \IN} \subset \cF$
of subsets $\Lambda_N$ whose cardinality 
does not exceed $N$, 
such that 
\be\label{eq:BestNgpcL2}
\left\| P - \sum_{\nu\in \Lambda_N} P_\nu L_\nu \right\|_{L^2(U,\rho;V)}^2
= 
\sum_{\nu\in{\mathcal F}\setminus\Lambda_N}\|P_\nu\|_V^2
\le 
C N^{-2\sigma},
\ee
where the convergence rate $ \sigma = 1/p-1/2>1/2 $ 
and where the constant 
$C = \left\| (\|P_\nu\|_V)_{\nu\in \cF} \right\|^2_{\ell^p(\cF)}$ 
is bounded independently of $N$.
\end{proposition}
The best $N$-term approximations
\be\label{eq:BestNtAppr}
P_{\Lambda_N}  := \sum_{\nu\in \Lambda_N} P_\nu L_\nu
\ee
in Proposition \ref{prop:Stechkin} indicate that
sampling the parametric forward map with evaluation
of $N$ solutions $P_\nu$, $\nu\in \Lambda_N$ 
of the parametric, elliptic problem with 
accuracy $N^{-\sigma}$ is possible; since $\sigma>1/2$, 
this is superior to what can be expected from $N$ MC samples.
There are, however, two obstacles which obstruct the practicality
of this idea: first, the proof of Proposition \ref{prop:Stechkin} is  
nonconstructive, and does not provide concrete choices for 
the sets $\Lambda_N$ of ``active'' gpc coefficients which realize 
\eref{eq:BestNgpcL2} and, second, even if $\Lambda_N$ were
available, the ``coefficients'' $P_\nu\in V$ can not be obtained
exactly, in general, but must be approximated for example
by a Finite Element discretization in $D$.

As $P\in L^2(U,\rho; V)$, 
we consider the variational form ``in the mean'' 
of \eref{eq:prob} as
\be
\fl\int_U\int_DK(x,u)\nabla P(x,u)\cdot\nabla Q(x,u)dxd\rho(u)
=
\int_U\int_Df(x)Q(x,u)dxd\rho(u),
\label{eq:parametricweak}
\ee
for all $Q\in L^2(U,\rho;V)$. 
For each set $\Lambda_N\subset{\mathcal F}$ 
of cardinality not {exceeding} $N$ 
that satisfies Proposition \ref{prop:Stechkin}, 
and each   vector ${\mathcal L}=(l_\nu)_{\nu\in\Lambda_N}$ of nonnegative 
integers, we define finite dimensional approximation spaces as
\be
X_{N,{\mathcal L}} 
= 
\{P_{N,{\mathcal L}}=\sum_{\nu\in\Lambda_N}P_{\nu,{\mathcal L}}(x)L_\nu(u);
\ 
P_{\nu,{\mathcal L}}\in V^{l_\nu}\}
\;.
\label{eq:XNL}
\ee
Evidently, $X_{N,{\mathcal L}} \subset L^2(U,\rho;V)$ 
is a finite-dimensional (hence closed) 
subspace for any $N$ and any selection ${\mathcal L}$ 
of the discretization levels.

The total number of degrees of freedom, 
$N_{dof} = {\rm dim} (X_{N,{\mathcal L}})$,
necessary for the sparse representation of the 
parametric forward map is given by
\begin{equation}\label{eq:DeftotNdof}
\eN = O\left(\sum_{\nu\in\Lambda_N}2^{dl_\nu}\right)
\quad\mbox{as} \quad N, \;l_\nu \rightarrow \infty
\;. 
\end{equation}
The 
{\em stochastic, sparse tensor Galerkin approximation}
of the parametric forward problem \eref{eq:prob}, 
based on the 
index sets $\Lambda_N\subset \cF$, and $\cL = \{ l_\nu : \nu\in \Lambda_N \}$,
reads:
find $P_{N,{\mathcal L}}\in X_{N,{\mathcal L}}$ 
such that for all $Q_{N,{\mathcal L}} \in   X_{N,{\mathcal L}}$ holds
\be \label{eq:parameterfe}
\begin{array}{rcl}
\ds
b(P_{N,{\mathcal L}},Q_{N,{\mathcal L}})
& := & \ds
\int_U\int_DK(x,u)\nabla P_{N,{\mathcal L}}\cdot\nabla Q_{N,{\mathcal L}}dxd\rho(u) 
\\
& = & \ds
\int_U\int_D f(x)Q_{N,{\mathcal L}}(x,u)dxd\rho(u)
\;.
\end{array}
\ee
The coercivity of the bilinear form $b(\cdot,\cdot)$ ensures the existence and uniqueness
of $P_{N,{\mathcal L}}$ as well as their quasioptimality in $L^2(U,\rho;V)$:
by Cea's lemma, for a constant $C>0$ which is independent of $\Lambda$ and of $\cL$,
\[
\|P-P_{N,{\mathcal L}}\|_{L^2(U,\rho;V)}
\le 
C\inf_{Q_{\nu,{\mathcal L}}\in V^{l_\nu}}
\|P-\sum_{\nu\in\Lambda}Q_{\nu,{\mathcal L}}L_\nu\|_{L^2(U,\rho;V)}
\;.
\]
We obtain the following error bound which consists
of the error in the best $N$-term truncation
for the gpc expansion and 
of the Finite Element approximation error for the 
``active'' gpc coefficients.
\be
\|P-P_{N,{\mathcal L}}\|_{L^2(U,\rho;V)}^2 
\le 
C(N^{-2\sigma}
+ 
\sum_{\nu\in\Lambda_N}\inf_{Q_{\nu,{\mathcal L}}\in V^{l_\nu}}\|P_{\nu}-Q_{\nu,{\mathcal L}}\|_V^2)
\;.
\label{eq:est1}
\ee
%

Let us indicate sufficient conditions that ensure 
Assumptions \ref{ass:gpclp}, \ref{ass:gpcplxity}.
The first condition is quantitative 
decay rate of the coefficient functions $\psi_j$ 
in the parametric representation \eref{eq:a}
of the random input.

To obtain convergence rates 
for the FE-discretization in the domain $D$,
i.e. of the term 
$\|P_{\nu}-Q_{\nu,{\mathcal L}}\|_V $ in \eref{eq:est1}, 
it is also necessary
to ensure spatial regularity of the solution $P(x,u)$ of the parametric
problem \eref{eq:prob}. 
To this end we employ Assumptions \ref{assump:sumpsi}(iii). 
We remark that Assumptions \ref{assump:sumpsi}(iii) and (iv) 
imply Assumption \ref{assump:sumpsi}(ii) with $q = s-1 > 0$.
Under these assumptions, the following proposition holds.
\begin{proposition} \label{prop:H2reg}
{
Under Assumptions \ref{assump:sumpsi}(i),(iii) and (iv) and
if, moreover, the domain $D$ is convex and $f\in L^2(D)$,
}
the solution $P(\cdot,u)$ of the 
{
parametric, deterministic} 
problem \eref{eq:prob} belongs to the space $L^2(U,\rho; W)$.
\end{proposition}
From estimate \eref{eq:est1}, we get with Proposition \ref{prop:H2reg}
and standard approximation properties of continuous, piecewise linear FEM
the error bound
\be
\|P-P_{N,{\mathcal L}}\|_{L^2(U,\rho;V)}^2
\le
C(N^{-2\sigma} + \sum_{\nu\in\Lambda_N}2^{-2l_\nu}\|P_{\nu}\|_{H^2(D)}^2)
\;.
\label{eq:est2}
\ee
In order to obtain an error bound in terms of $N_{dof}$ 
defined in \eref{eq:DeftotNdof} which is uniform in terms of $N$,
we select, for $\nu\in\Lambda_N$ 
the discretization levels $l_\nu$ of the active gpc coefficient 
$P_\nu$ so that both terms in the upper bound \eref{eq:est2} 
are of equal order of magnitude. 
This constrained optimization problem was solved, for example, 
in \cite{CDS1}, under the assumption that 
$(\| P_\nu \|_{H^2(D)})_\nu\in \ell^p(\cF)$.

In recent years, several algorithms have appeared
or are under current development 
which satisfy Assumption \ref{ass:gpcplxity} with 
various exponents $\alpha \geq 1$ and $\beta \geq 0$.
See for example the references 
\cite{BAS,ScCJG,CJGdiss,b:09,CCDS11}

\ack VHH is supported by a start up grant from Nanyang Technological University, CS acknowledges partial support by the
Swiss National Science Foundation
and by the
European Research Council under
grant ERC AdG 247277 - STAHDPDE and AMS is grateful to EPSRC (UK) and ERC for financial support.
The authors are also particularly grateful to Daniel Gruhlke, who shared
his ideas concerning multilevel MCMC methods for conditioned diffusion
processes during a visit to AS at Warwick in September 2011, and to
Rob Scheichl who found an error in an earlier pre-print of the multilevel
MCMC material, leading us to consider the improved analysis contained
herein.

\section*{References}
\bibliography{sparseMCMC}
\bibliographystyle{plain}
\end{document}